\documentclass[aop,MSNbibl,seceqn,rotating,citesort,dvips]{arximspdf}

\doi{10.1214/11-AOP672}
\volume{40}
\issue{5}
\pubyear{2012}
\firstpage{2168}
\lastpage{2196}

\makeatletter
\let\epsilon\varepsilon

\newtheorem{cor}[thm]{Corollary}% [section]
\newtheorem{lem}[thm]{Lemma}%[section]
\newtheorem{prop}[thm]{Proposition}%[chapter]
\newproclaim{rem}[thm]{Remark}

\newcommand{\eqref}[1]{(\ref{#1})}
\newcommand{\rr}{\varrho}
\newcommand{\QQ}{\mathbb{Q}}
\newcommand{\EE}{\mathbb{E}}
\newcommand{\HH}{\mathbb{H}}
\newcommand{\RR}{\mathbb{R}}
\newcommand{\PP}{\mathbb{P}}
\newcommand{\MMM}{\mathbb{M}}
\newcommand{\OO}{\mathbb{O}}
\newcommand{\VV}{\mathbb{V}}

\newcommand{\mcC}{\mathcal{C}}
\newcommand{\mcF}{\mathcal{F}}
\newcommand{\mcU}{\mathcal{U}}
\newcommand{\mcH}{\mathcal{H}}
\newcommand{\mcV}{\mathcal{V}}
\newcommand{\mcB}{\mathcal{B}}
\newcommand{\mcN}{\mathcal{N}}
\newcommand{\mcD}{\mathcal{D}}
\newcommand{\mcO}{\mathcal{O}}
\newcommand{\mcE}{\mathcal{E}}
\newcommand{\mcG}{\mathcal{G}}
\newcommand{\mcL}{\mathcal{L}}
\newcommand{\mcA}{\mathcal{A}}

\newcommand{\lar}{\leftarrow}

\newcommand{\be}{\mathbf{e}}
\newcommand{\bfT}{\mathbf{T}}

\newcommand{\Vol}{\mbox{\textsc{Vol}}}
\newcommand{\Vvol}{\mbox{\fontsize{8.36}{9}\selectfont\textsc{Vol}}}

\newcommand{\st}{such that }
\newcommand{\varep}{\varepsilon}
\newcommand{\ep}{\epsilon}

\newcommand{\la}{\lambda}
\newcommand{\al}{\alpha}

\newcommand{\wrt}{with respect to }
\newcommand{\wlg}{without loss of generality }
\newcommand{\bs}{\setminus}
\newcommand{\ra}{\rightarrow}
\newcommand{\p}{\partial}

\def\moins{\setminus}
\def\1{^{-1}}

\makeatother

\begin{document}
\begin{frontmatter}

\title{Nonexplosion criteria for relativistic diffusions\thanksref{TITL}}
\runtitle{Nonexplosion criteria for relativistic diffusions}
\thankstext{TITL}{Supported in part by the Agence Nationale de la
Recherche, ANR-09-BLANC-0364-01, and the EPSRC Grant EP/E01772X/1.}

\begin{aug}
\author[A]{\fnms{Isma\"{e}l} \snm{Bailleul}\corref{}\ead[label=e1]{i.bailleul@statslab.cam.ac.uk}\ead[url,label=u1]{http://www.statslab.cam.ac.uk/\textasciitilde ismael/}}
and
\author[B]{\fnms{Jacques} \snm{Franchi}\ead[label=e2]{jacques.franchi@math.unistra.fr}\ead[url,label=u2]{www-irma.u-strasbg.fr/\textasciitilde franchi}}
\runauthor{I. Bailleul and J. Franchi}
\affiliation{Centre for Mathematical Sciences and Universit\'{e} de
Strasbourg et CNRS}
\address[A]{Statistical Laboratory\\
Centre for Mathematical Sciences\\
Wilberforce Road\\
Cambridge CB3 0WB\\
United Kingdom\\
\printead{e1}\\
\printead{u1}} %adresu isvedimo komanda gale!
\address[B]{Universit\'{e} de Strasbourg et CNRS\\
I.R.M.A.\\
7 rue Ren\'{e} Descartes\\
67084 Strasbourg cedex\\
France\\
\printead{e2}\\
\printead{u2}}
\end{aug}

% HISTORY:
\received{\smonth{8} \syear{2010}}
\revised{\smonth{3} \syear{2011}}

% ABSTRACT
%
\begin{abstract}
Some general Lorentz covariant operators, associated to the so-called
$\Theta$ (or $\Xi$)-relativistic diffusions and making sense in any
Lorentzian manifold, have been introduced by Franchi and Le Jan
[\textit{Comm. Pure Appl. Math.} \textbf{60} (2007) 187--251],
Franchi and Le Jan [Curvature diffusions in general relativity (2010).
Unpublished manuscript]. Only a few examples have been studied so far.
We provide in this work some nonexplosion criteria for these
diffusions, which can be used in generic cases.
\end{abstract}

% KEYWORDS
%

\end{frontmatter}

%s1 ###
\section{Introduction}
\label{SectionIntroduction}
%---------------------%

It is well known that the metric completeness of a Riemannian manifold
does not prevent Brownian motion  from exploding within a finite time with positive
probability. The situation is now well understood, in particular,
thanks to the works of Yau~\cite{YauHeatKernel}, Grigor'yan
\cite{GrigoryanCompleteness}, Takeda~\cite{Takeda1,Takeda2}
and very recently Hsu and Qin~\cite{HsuQin}, to cite but a few names.
Different lines of approach have been used. Yau and Grigor'yan treated
the analytic counterpart of the completeness problem and investigated
the well-posedness of the parabolic Cauchy problem; the former using
local information on the geometry under the form of curvature bounds;
the latter using a global information under the form of an upper bound
for the volume of large balls. Takeda used a purely probabilistic
method due to Lyons and Zheng in~\cite{LyonsZheng}, based on
reversibility. This approach was recently improved by Hsu and Qin in
\cite{HsuQin}. Hsu used stochastic analysis in~\cite{HsuBook}, Theorem~3.5.1, to control the radial process, by estimating the Laplacian of
the distance function to a~fixed point in terms of curvature bounds.
All these results are tied down to the metric framework provided by a
complete Riemannian manifold.

A natural analog of Brownian motion  in a Lorentzian setting was first
introduced by Dudley~\cite{Dudley1} in the special relativistic case,
and extended to the general relativistic framework by Franchi and Le
Jan \cite {FranchiLeJan}. It belongs to a larger class of relativistic
processes introduced in~\cite{BailleulIHP} and
\cite{FranchiLeJanCurvature}, defined in purely geometric terms and
collectively refered to as \textit{relativistic diffusions}. Their
trajectories represent the random motion in spacetime of a~small
massive particle, and make sense only in the unit tangent bundle or in
the orthonormal frame bundle. Only a few examples have been studied in
detail up to now: in Minkowski spacetime (the framework of special
relativity)~\cite{Dudley1,BailleulPoisson,BailleulRaugi}, in
Robertson--Walker spacetimes (models of universe with a big-bang)
\cite{JurgenPhD}, G\"{o}del spacetime (a causally paradoxical universe)
\cite{FranchiGodel} and Schwarzschild spacetime (a model for an
isolated star or a black hole)~\cite{FranchiLeJan}.

Apart from the works~\cite{BailleulIHP} and \cite
{FranchiLeJanCurvature}, no general study of these intrinsic random
processes was done. As a first step towards a better understanding of
these processes and their interplay with the geometry of the ambient
spacetime, we provide in this work some nonexplosion criteria for some
generic classes of Lorentz manifolds. In addition to being a natural
question, the completeness issue is strongly related to important
questions in general relativity. Indeed, dating back to Penrose and
Hawking's incompleteness theorems, the appearance of singularities in
Einstein's theory of gravitation has been recognized as unavoidable
under quite natural assumptions. Although there is no agreement on what
should be called a singularity of a spacetime, the existence of
incomplete geodesics has been widely used as an indicator of such a
singular feature. In so far as the random dynamics considered in this
work (Section~\ref{SectionRanDyn}) can be seen as intrinsic
perturbations of the geodesic flow, their completeness/incompleteness
is a distinguishing feature of a spacetime. We refer the reader to
\cite{BailleulJMP} for a first approach of stochastic incompleteness.

The paths of the random processes we shall consider are (almost-)all
$\mcC^1$ paths parametrized by their (proper time) arc length. What
could possibly make them explode? In a complete Riemannian manifold,
any such path would have to be at time $s$ in a closed ball of radius
$s$ with center its starting point, so it cannot explode. There are two
problems with the Lorentzian setting: a Lorentzian manifold has no
metric or finite distance function associated with its structure, and
the set of unit tangent vectors at any point is noncompact. As a
result, even in Minkowski spacetime, one can construct exploding paths
with finite (proper time) arc length.

To start our investigations, we shall take advantage in Section \ref
{SectionLiapounov} of the bundle structure of the state space of the
process, to exhibit a one-dimension\-al sub-process whose control is
possible in the class of globally hyperbolic spacetimes. This
structure, indeed, allows us to define some Lyapounov function, and
leads to a nonexplosion criterion by using a simple and well-known
observation due to Khasminsky.

With a metric missing, the completeness notion used in a crucial way in
the Riemannian setting becomes unavailable. Busemann, Hawking and Ellis
and Schmidt, Beem and Ehrlich proposed different notions in
replacement. Schmidt's\vadjust{\goodbreak} idea is to give a Riemannian structure to the
orthonormal frame bundle. We consider Schmidt b-completeness notion
in Section~\ref{SectionbCompleteness}, showing how it leads to a
stochastic completeness result for some of the relativistic diffusions.

This result can be significantly improved by adapting Takeda's
strategy~\cite{Takeda2}, as improved by Hsu and Qin~\cite{HsuQin}, to the
Lorentzian setting. This is, however, far from being straightforward,
since we are working in a nonsymmetric, nonelliptic setting, where the
main ingredients of Takeda's method (use of symmetry and reflected
Brownian motion on the boundary of large Riemannian balls) have no
obvious Lorentzian counterpart. To overcome this difficulty, we use in
Section~\ref{SectionVolumeGrowth} a sub-Riemannian structure well
adapted to our setting, and which will somehow play for us the role of
the nonexisting Lorentzian distance.

%-------------------------------%
%s2 ###
\section{Relativistic diffusions}
\label{SectionThetaDiffusions}
%-------------------------------%

%%-----------------------------------%
%s2.1 ###
\subsection{Basic geometrical setting}
\label{SectionGeometricalSetting}
%%-----------------------------------%

Recall Minkowski space is the product $ \RR^{1,d}\equiv\RR\times\RR
^d$ equipped with the metric
\[
g_{\mathrm{M}}(q,q) := t^2 - |x^1 |^2 - \cdots- |x^d |^2 \qquad  \mbox{for
any } q=(t,x)\in\RR^{1,d},
\]
where $(t,x^1,\dots,x^d)$ denote the coordinates of $ q$ in the
canonical basis $ \{\ep_0,\break\ep_1,\dots,\ep_d \}$ of $\RR^{1,d}$.

Let $(\MMM,g)$ be a smooth $(1+d)$-dimensional Lorentzian manifold
(with $ d\ge2$), which we shall always suppose to be oriented and
time-oriented. (We refer the reader to the books of Hawking--Ellis
\cite{HawkingEllis} and O'Neill~\cite{Oneill} for the basics on
Lorentzian geometry.) Given any point $m\in\MMM$, it is usual to
consider an orthonormal basis $\{\be_0,\ldots,\be_d\}$ of the tangent
space $ T_m\MMM$ as an isometry $ \be$ from $ (\RR^{1,d},g_{\mathrm{M}} )$ to $ (T_m\MMM,g_m )$; so, strictly speaking, $\be_i = \be
(\ep_i)$. The orthonormal frame bundle of $\MMM$ is just the collection
\[
\OO\MMM= \{\Phi= (m,\be) | m\in\MMM, \be \mbox{ an
orthonormal basis of }(T_m\MMM, g_m) \}.
\]
We shall write $\OO\mcU= \{\Phi= (m,\be) | m\in\mcU, \be\mbox{ an orthonormal basis of }T_m\MMM\}$ for any subset $\mcU$ of $\MMM
$. For a small enough $\mcU$ and a chart $x \dvtx \mcU\ra\RR^{1+d}$ on
it, we shall write $ \be_j = e_j^k\,\partial_{x^k} $ for each vector $
\be_j $ of a frame $ \be$; this decomposition provides local
coordinates $(x^i,e_j^k)$ on $\OO\mcU$.\vspace*{2pt}

Each fiber $\OO_m\MMM$ is modeled on the noncompact orthogonal group
$O(1,d)$, which has four connected components. We shall be interested
in dynamics leaving these components globally fixed. We choose to
consider only one of them, specified by the requirement that $\be_0$
should be future-oriented and that the orientation of $\be$ should be
direct. We shall still denote the resulting frame bundle by $\OO\MMM
$, as there will be no risk of confusion. The Lorentz--M\"{o}bius group
$\mathit{SO}_0(1,d)$, that is, the connected component of the unit in $O(1,d)$,
acts properly on $\OO\MMM$. This natural action induces the canonical
vertical vector fields $(V_{ij})_{0\leq i<j\leq d}$. The subgroup of
elements in $\mathit{SO}_0(1,d)$ that\vadjust{\goodbreak} fix $ \ep_0 $ can be identified with the
rotation group $\mathit{SO}(d)$, and generates the vector fields
$(V_{ij})_{1\leq i<j\leq d}$. To shorten notations we shall write $V_j$
for $V_{0j}$; it generates boosts, that is, hyperbolic rotations in
each fiber, and reads, in the above local coordinates,
%
%e2.1 ###
\begin{equation} \label{f.Vj}
V_j = e_j^k \frac{\partial}{\partial e_0^k} + e_0^k \frac{\partial
}{\partial e_j^k} .
\end{equation}

Throughout this work, $T\MMM$ and $\OO\MMM$ will be endowed with the
Levi--Civita connection, inherited from the Lorentzian pseudo-metric
$g$. Last, we denote by $ H_0 $ the vector field generating the
geodesic flow on $\OO\MMM$. Denoting by $\Gamma^\ell_{kj}$ the
Christoffel coefficients, we have, in the above local chart on $\OO
\MMM$,
%
%e2.2 ###
\begin{equation} \label{f.H0}
H_0 = e_0^k \,\partial_{x^k} - e_0^k e_i^j \Gamma^\ell_{kj} \frac
{\partial}{\partial e_i^\ell} .
\end{equation}

We shall denote by $T^1\MMM$ the future-oriented unit tangent bundle
over $\MMM$, with generic element $(m,{\dot m})$. In Minkowski
spacetime $\RR^{1,d}$, it is the product of $\RR^{1,d}$ by the
hyperboloid $\HH= \{q=(t,x)\in\RR^{1,d} ; g(q,q)=1, t>0 \}$. The
bundle $T^1\MMM$ is locally modeled on that product. (Consult \cite
{HawkingEllis} or~\cite{Oneill} for some background.) Denote by $ \pi
_1 $ the projection $(m,\be)\mapsto(m,\be_0\equiv\dot m)$ from $\OO
\MMM$ to $T^1\MMM$, and by $ \pi_0 $ the canonical projection $\OO
\MMM\ra\MMM$.

%%--------------------------------------%%
%s2.2 ###
\subsection{Relativistic random dynamics}
\label{SectionRanDyn}
%%--------------------------------------%%

Relativistic diffusions model the random motion in spacetime of a small
massive particle parametrized by its proper time, providing random
timelike paths$ $; so, properly speaking, their mathematical
counterpart are random trajectories $(m_s,\dot m_s)$ in the future unit
bundle $T^1\MMM$ subject to the condition $\frac{d}{ds}m_s = \dot m_s
$. Yet it happens to be more convenient to define random dynamics in
the orthonormal frame bundle $\OO\MMM$ as it bears more structure
than $T^1\MMM$; these diffusions on $\OO\MMM$ are constructed so as
to have a projection on $T^1\MMM$ which is itself a diffusion. Such a
construction is remniscent of Malliavin--Eells--Elworthy's construction
of Brownian motion  on a Riemannian manifold as the projection of a~diffusion on the
orthonormal frame bundle.

%%%-------------------------------%%%
%s2.2.1 ###
\subsubsection{Dynamics in $\OO\MMM$}
%%%-------------------------------%%%

Given any smooth nonnegative function $\Theta\dvtx\break T^1\MMM\ra\RR_+$,
identified to a $\mathit{SO}(d)$-invariant function on $\OO\MMM$ by setting $
\Theta(\Phi):=\Theta(\pi_1(\Phi) )$, consider the following
Stratonovich differential equation on $\OO\MMM$:
%
%e2.3 ###
\begin{eqnarray}
\label{RelativisticDynamics}
\circ \,d\Phi_s &=& H_0(\Phi_s) \,ds + \frac{1}{4} \sum_{1\le
j\le d} V_j\Theta(\Phi_s) V_j(\Phi_s) \,ds \nonumber\\[-8pt]\\[-8pt]
&&{}+ \sqrt{\Theta(\Phi_s)}\sum_{1\le j\le d} V_j(\Phi_s) \circ \,dw^j_s,\nonumber
\end{eqnarray}
where $ w$ is a $d$-dimensional Brownian motion and where we understand
a~vector field as a first-order differential operator.  This equation
has a unique maximal strong solution, defined up to its explosion time
$\zeta$.\vadjust{\goodbreak}

It is clear on this equation that the $(\be_1,\ldots,\be_d)$-part of
$\Phi_s$ is irrelevant in defining the dynamics of $ (m_s,\be_0(s) )$
since $\Theta(\Phi)$ depends only on $\pi_1(\Phi) $; this is the
reason why this diffusion on $\OO\MMM$ projects down in $T^1\MMM$
onto a~diffusion. Consult~\cite{FranchiLeJan}, Theorem $1$, \cite
{FranchiLeJanCurvature}, Theorem $3.2.1$ or~\cite{BailleulIHP},
Section 3.2, for the details. The diffusion in $\OO\MMM$ has generator
%
%e2.4 ###
\begin{equation} \label{f.Gen}
\mcG_{\Theta} = H_0 + \frac{1}{2} \sum_{1\le j\le d} V_j
(\Theta V_j ) .
\end{equation}

We shall generically call these relativistic dynamics \textit{$ \Theta$-diffusions} (the $\Xi$-diffusions of \cite
{FranchiLeJanCurvature}). These diffusions are \textit{covariant}, in the sense that any isometry of $(\MMM,g)$ maps a $
\Theta$-diffusion to a $ \Theta$-diffusion (with the same $ \Theta$:
the law is preserved, up to the starting point), and admit the
Liouville measure as an invariant measure. The $ \pi_0$-projections
(on the base manifold $\MMM$) of their trajectories are almost-surely
$C^1$ paths.  A $ \Theta$-diffusion $(\Phi_s)_{0\le s<\zeta}$
solving equation~(\ref{RelativisticDynamics}) is parametrized by
proper time $ s\ge0 $.  The particular case $ \Theta=0 $ gives back
the deterministic geodesic flow, and the case of a nonnull constant
$\Theta$ gives back the relativistic diffusion as defined first in
\cite{FranchiLeJan}, which we shall call the \textit{basic
relativistic diffusion}. It is described in simple terms in Minkowski
spacetime. Although the metric $g_M$ is \mbox{nondefinite} positive, its
restriction to any tangent space of the half sphere $\HH$ of unit
tangent vectors is definite negative; this turns $\HH$ into
a~Riemannian manifold with constant negative curvature. Dudley's
diffusion $ (m_s,\be_s ) = (m_s,(\be_0(s),\dots,\be_d(s)) )$, which
is the basic relativistic diffusion in Minkowski spacetime, corresponds
to taking $m_s = m_0+\int_0^s \be_0(r) \,dr$, and for the velocity $\be
_0(r)$ a Brownian motion  on $\HH$. The remainder $\be_1(r),\dots, \be_d(r)$ of
the basis is obtained by paralell transport of $\be_1(0),\dots, \be
_d(0)$ along the Brownian path $(\be_0(u))_{0\leq u\leq r}$.

The following elementary lemma, proved in~\cite{BailleulJMP}, Section
2.2, gives an intuitive picture of the $\Theta$-diffusions, for
$\Theta$ depending only on $m\in\MMM$.

\begin{lem}\label{LemmaLifting}
Let $\gamma\dvtx [0,T]\ra\MMM$ be a $\mcC^2$ timelike path parametrized
by its proper time, and $\Gamma_0\in\OO\MMM$ \st$\pi_1(\Gamma_0)
= (\gamma(0),\dot\gamma(0) )\in T^1\MMM$. Then there exists a
unique $\mcC^2$ path $ (\Psi_s )_{0\leq s\leq T}$ in $\OO\MMM$, and
some unique $\mcC^1$ real-valued controls $h^1,\dots,h^d$ defined on
$[0,T]$, \st$\Psi_0 = \Gamma_0, \pi_1(\Psi_s) = (\gamma(s),\dot
\gamma(s) )$ and
\[
\dot\Psi_s = H_0(\Psi_s) + \sum_{j=1}^d V_j(\Psi_s) h^j(s).
\]
\end{lem}

So the $\Theta$-diffusion is obtained in that case by replacing
the deterministic controls of a typical $\mcC^2$ timelike path by
Brownian controls with position dependent variance
$\Theta(m_s)$.

On a manifold with nonpositive scalar curvature $ R $, taking $ \Theta
(\Phi) = -\rr^2 R$ (for a nonnull constant $\rr$), one gets a
dynamic which can be truly random only in nonempty parts of spacetime$
$; it was called \textit{$ R$-diffusion} in~\cite
{FranchiLeJanCurvature}. Denote by $\bfT$ the energy-momentum tensor
of the spacetime. Taking $ \Theta(\Phi) = \rr^2 \bfT(\be_0,\be
_0)$, we get what was named the \textit{energy diffusion} in
\cite{FranchiLeJanCurvature}. See~\cite{BailleulIHP} for more general
models of diffusions.

%%%-------------------------------%%%
%s2.2.2 ###
\subsubsection{Dynamics in $T^1\MMM$}
%%%-------------------------------%%%

Denote by $\nabla^v$ the gradient on $T^1_m\MMM$, identified with the
hyperbolic space $\HH^d$ by means of the metric $ g_m $, and by $\mcL
_0$ the vector field generating the geodesic flow on $T^1\MMM$. Note
that $ T\pi_1(H_0) = \mcL_0$ and $T\pi_1(V_j) =: \nabla^v_j = e^k_j
\partial_{{\dot m}^k}$ (with Einstein summation convention). The
projection on $T^1\MMM$ of the $\OO\MMM$-valued diffusion has the
following $\mathit{SO}(d)$-invariant generator:
\[
\mcL_{\Theta} = \mcL_0 + \tfrac{1}{2} \nabla^v (\Theta\nabla
^v ) .
\]
For a constant $\Theta$ the operator $\mcL_{\Theta}$ has the
following expression in the local coordinates introduced in Section
\ref{SectionGeometricalSetting}:
\begin{eqnarray*}
\mcL_0 + \frac{\Theta}{2} \Delta^v &=& {\dot m}^k\frac{\partial
}{\partial{m}^k} + \biggl( \frac{d}{2} \Theta{\dot m}^k - {\dot
m}^i{\dot m}^j \Gamma_{ij}^k(m) \biggr) \frac{\partial}{\partial{\dot
m}^k}\\
&&{} + \frac{\Theta}{2} \bigl({\dot m}^k{\dot m}^\ell- g^{k\ell}(m) \bigr)
\frac{\partial^2}{\partial{\dot m}^k \,\partial{\dot m}^\ell} ,
\end{eqnarray*}
where $\Delta^v$ denotes the vertical Laplacian. We have, for a
generic $\Theta$,
%
%e2.5 ###
\begin{equation}
\label{GeneratorCoordinates}
\mcL_{\Theta} = \mcL_0 + \frac{\Theta}{2} \Delta^v + \frac
{1}{2} \bigl({\dot m}^k{\dot m}^\ell-g^{k\ell}(m) \bigr) \frac{\partial
\Theta}{\partial{\dot m}^k} \frac{\partial}{\partial{\dot m}^\ell
} .
\end{equation}

The purpose of this work is to provide some conditions under which the
$\Theta$-diffusions have almost-surely an infinite lifetime $ \zeta$. In so far
as we are mainly interested in the $T^1\MMM$-valued $\Theta
$-diffusions as models of physical phenomena, while we shall mainly be
working with $\OO\MMM$-valued diffusions, it is reassuring to have
the following fact, which essentially means that the possible explosion
of $(\Phi_s)_{0\le s<\zeta}$ is never due to its $(\be_1,\ldots,\be_d)$-part.

\begin{prop}
The $\Theta$-diffusion on $\OO\MMM$ and its $T^1\MMM$-projection
have the same lifetime.
\end{prop}

\begin{pf}
Write $\Phi_s= (m_s ; ({\dot m}_s, e_1(s),\ldots,e_d(s) ) )\in\OO
\MMM$ and $ \phi_s := \pi_1(\Phi_s)= (m_s,{\dot m}_s)\in T^1\MMM$.
Using the local coordinates $(x^k,e_j^\ell)_{0\leq k,\ell\leq d ;
1\leq j\leq d} $, equation~\eqref{RelativisticDynamics} defining the
$\Theta$-diffusion reads
\begin{eqnarray*}
d{\dot m}_s^k &=& dM^k_s -\Gamma_{i\ell}^k(m_s) {\dot m}_s^i {\dot
m}_s^\ell \,ds + \frac{d}{2} \Theta(\phi_s) {\dot m}_s^k \,ds \\
&&{}+
\frac{1}{2} \bigl({\dot m}^k_s{\dot m}^\ell_s - g^{k\ell}(m_s) \bigr)
\frac{\partial\Theta}{\partial{\dot m}^\ell}(\phi_s) \,ds , \\
de^k_j(s) &=& \sqrt{\Theta(\phi_s)} {\dot m}_s^k \,dw^{j}_s - \Gamma
_{i\ell}^k(m_s) e^\ell_j(s) {\dot m}_s^i \,ds \\
&&{}+ \frac{1}{2}
\Theta(\phi_s) e^k_j(s) \,ds + \frac{1}{2} V_j \Theta(\phi_s)
{\dot m}_s^k \,ds ,
\end{eqnarray*}
with the martingale term $ dM^k_s := \sqrt{\Theta(\phi_s)} e^k_j(s)
\,dw^{j}_s$. (See Section 3.2 of~\cite{FranchiLeJanCurvature} for the
computation of the It\^{o} correction.) Setting $e_0={\dot m}$ and
$\eta^{in} :=\eta_i^n := \mathbf{1}_{i=n=0}-\mathbf{1}_{1\le i=n\leq d} $,
and noticing that the matrix $ (\eta^{in} e_n^k g_{k\ell} )_{0\leq
i,\ell\leq d}$ is the inverse of the matrix $ (e^{i}_{\ell} )_{0\leq
i,\ell\leq d} $, it follows from the above system that we have, for
all $0\le k\le d , 1\le j\le d$,
\begin{eqnarray*}
de^k_j(s) &= &\dot m_s^k \eta^{n}_j e_n^q(s) g_{q\ell}( m_s)\, dM^{\ell
}_s -\Gamma_{i\ell}^k( m_s) e^\ell_j(s) \dot m_s^i \,ds + \frac{1}{2} \Theta(\phi_s) e^k_j(s) \,ds\\
&&{} + \frac{1}{2} V_j\Theta
(\phi_s) \dot m_s^k \,ds \\
&=& - e^\ell_j(s) \Gamma_{i\ell}^k( m_s) \dot m_s^i \,ds + \frac{1}{2} e^k_j(s) \Theta(\phi_s) \,ds + \frac{1}{2} V_j\Theta
(\phi_s) \dot m_s^k \,ds \\
&&{} - e_j^q(s) \dot m_s^k g_{q\ell}( m_s) \biggl[d\dot m_s^\ell+\Gamma
_{ip}^\ell( m_s) \dot m_s^i \dot m_s^p \,ds - \frac{d}{2} \Theta
(\phi_s) \dot m_s^\ell \,ds\\
&&\hspace*{116pt}
{} - \frac{1}{2} [\dot m^p_s \dot
m^\ell_s - g^{p\ell}( m_s) ] \frac{\partial\Theta}{\partial\dot
m^p}(\phi_s) \,ds \biggr] .
\end{eqnarray*}
So the matrix $ (e^{k}_{j}(s) )_{0\leq s<\zeta}$ and the frame-valued
diffusion $(\Phi_s)_{0\leq s<\zeta} $ satisfy a linear stochastic
differential equation, conditionally on $(\phi_s)_{0\leq s<\zeta} $.
It is thus well defined up to the explosion time $ \zeta$ of the
$T^1\MMM$-valued $\Theta$-diffusion.
\end{pf}

 This point being clarified, we shall work freely in the sequel
with $\Theta$-diffusions on $\OO\MMM$.

%-----------------------------------------%
%s3 ###
\section{A first nonexplosion criterion}
\label{SectionLiapounov}
%-----------------------------------------%

We give in this section a simple nonexplosion criterion, well suited to
investigate the behavior of the $\Theta$-diffusions in the largely
used class of globally hyperbolic spacetimes. A Lyapounov function is
introduced for this purpose, and leads to a nonexplosion criterion of a
different nature than the typical Riemannian criteria mentioned in the
\hyperref[SectionIntroduction]{Introduction}.

The idea is roughly the following: if we can find a function $ f=f(\Phi
)$ which has compact level sets $\{f\leq\la\}$, and does not increase
along the trajectories, then the dynamics cannot explode. This was
noted first by Khasminsky in a~stochastic context; we state his
observation here for the relativistic diffusions.

\begin{lem}[(Khasminsky)] \label{lem.Khasm}
If there exists a nonnegative function $f$ on $\OO\MMM$ and a
positive constant $C$ \st$ \mcG_\Theta f\leq C f$, and $f$ goes to
infinity along any timelike path leaving any compact in a finite time,
then the $\Theta$-diffusion has almost-surely an infinite lifetime.\vadjust{\goodbreak}
\end{lem}

\begin{pf}
The condition $ \mcG_\Theta f\leq C f $ implies that the real-valued
process $ ( e^{-Cs}f(\Phi_s) )_{s<\zeta}$ is a nonnegative
supermartingale. Denote by $ \tau_n $ the (possibly infinite) exit
time from the level set $\{f\leq n\}$. By optional stopping, we have
\[
f(\Phi_0)\geq\EE[e^{-C \tau_n}f(\Phi_{\tau_n}) ] = n \EE[e^{-C
\tau_n} ].
\]
This implies that $\tau_n$ goes to infinity as $n$ goes to infinity;
as $ \zeta= \lim_{n\ra\infty} \tau_n $, this proves
Khasminsky's statement.
\end{pf}

As $\Theta$-diffusions have no a priori reason not to explode, such a
Lyapounov function will generally not exist. Yet, it is possible to
construct such a function in some classes of spacetimes of interest for
cosmology and theoretical physics. We give below two such examples. The
construction of the function~$f$ uses the same recipe in both cases:
if there exists an intrinsic distinguished future-directed timelike
$C^1$ vector field $ U\in T^1\MMM$, we can define
%
%e3.1 ###
\begin{equation}
\label{Definitionf}
f(\Phi) := g_m(U_m,\dot m);
\end{equation}
recall that $ \pi_1(\Phi) = (m,\dot m)\in T^1\MMM$. For this choice
of $f(\Phi)$, which is the hyperbolic angle between $ U$ and $\dot m
$, we have $ f\geq1$, and
%
%e3.2 ###
\begin{equation}
\label{FormulaH_0f}
H_0f(\Phi) = \nabla_{\dot m} ( g(U,\dot m) ) = g(\nabla_{\dot
m}U,\dot m) .
\end{equation}
The following lemma shows why $f$ is a good choice to apply
Khasminsky's criterion.

\begin{lem}
\label{LemmeSimple}
We have on $ \OO\MMM\dvtx\frac{1}{2} \sum_{j=1}^d
V_j(\Theta V_jf) = \frac{d}{2} \Theta f + \frac{1}{2} (f
\dot m^k-U^{k}) \frac{\partial\Theta}{\partial\dot m^k} $.
\end{lem}

\begin{pf}
Choose local coordinates for which $U=\partial_{x^0}$, so $ f(\Phi) =
\dot m^0 = e_0^0 $. Using~(\ref{f.Vj}), we have thus locally:
\[
V_jf =\biggl ( e^k_j \frac{\p}{\p e^k_0}+e^k_0 \frac{\p}{\p e^k_j} \biggr)
e^0_0 = e_j^0 ,\qquad  V_j^2 f = e^0_0 = f
\]
and
\[
\hspace*{-25pt}\sum_{j=1}^d (V_j\Theta)(V_jf) = \sum_{j=1}^d e_j^0 e_j^k \frac
{\partial\Theta}{\partial\dot m^k} = (\dot m^0 \dot m^k\!-\!g^{0k})
\frac{\partial\Theta}{\partial\dot m^k} = (f \dot m^k\!-\!U^{k}) \frac
{\partial\Theta}{\partial\dot m^k} .
\]
%
% Fix $m\in\MMM$. As each vector field $V_i$ leaves the fiber $T^1_m
%coordinates $(\rho,\sigma)\in\RR_+\times\mathbb{S}^{d-1}$ associated
%with $U\in T_m\MMM$. Then, $f(\be) = \textrm{ch}\rho$, and
%$$
%$$
\upqed\end{pf}

It follows from~\eqref{f.Gen} and~\eqref{FormulaH_0f} that
\[
\mcG_\Theta f = g(\nabla_{\dot m}U,\dot m) + \frac{d}{2} \Theta f +
\frac{1}{2} (f \dot m^k-U^{k}) \frac{\partial\Theta}{\partial
\dot m^k} \cdot
\]
Khasminsky's criterion will thus guarantee the nonexplosion of the
$\Theta$-diffusion provided $f$ explodes along exploding trajectories,
and there exists a positive constant $C$ \st
%
%e3.3 ###
\begin{equation}
\label{CriterionNonExplosion}
g(\nabla_{\dot m}U,\dot m) + \frac{1}{2} (f \dot m^k-U^{k})
\frac{\partial\Theta}{\partial\dot m^k} \leq \biggl(C- \frac
{d}{2} \Theta\biggr) g(U,\dot m).\vadjust{\goodbreak}
\end{equation}
In order to turn this criterion into an effective tool, we first
restrict ourselves to the following general class of spacetimes. This
inequality become s particularly simple when $\Theta$ depends only on
the base point $m\in\MMM$.

%%----------------------------------------%%
%s3.1 ###
\subsection{Globally hyperbolic spacetimes}
\label{SectionGloballyHyperbolic}
%%----------------------------------------%%

This class of cosmological models is characterized by the existence of
a global time function (i.e., a function $ \tau\dvtx \MMM\ra\RR$, with
timelike gradient) such that it has connected spacelike level sets $\{
\tau= t\}$ of $\tau$, and each integral curve of the vector field $
\nabla\tau$ meets each level set of $\tau$ in exactly one point.
Thus $\MMM$ is diffeomorphic to the product $I\times S$ of an interval
$I$ and a $d$-dimensional manifold $S$. Without loss of generality, we
can suppose the interval $I$ unbonded from above. With the example of
Minkowski spacetime in mind, we see that a given spacetime may have an
infinity of time functions; they are not intrinsically associated with
the geometry.

Yet, we can take for vector field $U$ in this setting the gradient of
the time function $\tau\dvtx m=(t,x)\in I\times S \mapsto t$, so
\[
f(\Phi) = g(U,\dot m) = \nabla_{\dot m} \tau= {\dot m}^0 = \dot t > 0.
\]
There is no hope, though, to prove inequality \eqref
{CriterionNonExplosion} without specifying further the model, as the
time function is not intrinsically defined. To proceed further, we
shall look at the sub-class of \textit{generalized warped
product spacetimes}, in which the time function is supplied by the
model and can be seen as an absolute time. These universes are globally
hyperbolic spacetimes $ \MMM= I\times S $ whose metric tensor has the form
%
%e3.4 ###
\begin{equation} \label{f.wp}
g_m(\dot m,\dot m) = a_m^2 |{\dot m}^0 |^2 - h_{m}({\dot m}^S,{\dot
m}^S) ,
\end{equation}
where $ {\dot m}^0$ is the image of $ \dot m\in T^1_m\MMM$ by the
differential of the first projection $ I\times S \ra I $ and $ {\dot
m}^S$ the image of $ \dot m $ by the differential of the second
projection $ I\times S \ra S$. Write  $m=(t,x)\in I \times S$. The
function $ a $ is a positive $\mcC^1$ function on $ \MMM$, assumed to
be bounded on any subset $I'\times S$ where $I'$ is bounded from above
and $ h_{m} $ is a positive-definite scalar product on $ T_xS$,
depending on $ m $ in a $ \mcC^1$ way. This class of spacetimes
contains all Robertson--Walker spacetimes (hence in particular de
Sitter and Einstein--de Sitter spacetimes and the universal covering
of the anti-de Sitter spacetime).

\begin{thm}
\label{ThmWarped}
Let $(\MMM,g)$ be a generalized warped product spacetime. If the function
\[
T^1\MMM\ni(m,{\dot m}) \longmapsto \nabla_{\dot m}\log a -
\frac{d}{4} \Theta(m,\dot m) - \frac{1}{4} \biggl(\dot m^k \frac
{\partial\Theta}{\partial\dot m^k } - \frac{1}{a^2(m) {\dot m}^0}
\frac{\partial\Theta}{\partial{\dot m}^0} \biggr)
\]
is bounded below, then the $\Theta$-diffusion almost-surely does not explode.
\end{thm}

\begin{pf}
$\bullet$ We first check that if the $\Theta$-diffusion has a finite
lifetime~$\zeta$ then $f(\Phi_s)$ explodes at time $\zeta^-$. To
that end, consider a timelike trajectory $ \gamma=(m_s,{\dot
m}_s)_{0\le s<T} = ((t_s,x_s),{\dot m}_s )_{0\le s<T}$ in\vadjust{\goodbreak} $T^1\MMM$,
defined on some semi-open interval $[0,T)$, and \st$\frac
{d}{ds}m_s={\dot m}_s $ and $ f(\gamma_s)=\dot t_s$ is bounded above
by some positive constant $C$. It follows that $ t_0\le t_s\le t_0+
CT$, and $ h_{m_s}(\dot x_s, \dot x_s) \leq C^2a^2_{m_s}$ is bounded
above by a constant since $a$ is bounded above on $(\inf
I,t_0+CT]\times S$. This entails that $(x_s)_{0\le s<T}$ cannot exit
a~bound\-ed region of $S$, and so that $ \gamma$ must be trapped in a
finite union of sets of the form $ \overline{J^+(m_0)}\cap\overline
{J^-(q_j)} $, for some $ q_j\in\MMM$. Such a union of sets is compact
in a hyperbolic spacetime (see, e.g.,~\cite{HawkingEllis},
Section 6.6), $\gamma$ is trapped in a~compact set. Would $\gamma$
explode, it would have a cluster point at which the strong causality
would fail, leading to a contradiction as globally hyperbolic
spacetimes are strongly causal (\cite{HawkingEllis}, Section 6.6).

$\bullet$ The condition of the theorem is a rephrasing of the local
condition~\eqref{CriterionNonExplosion}. To see that, let us work in a
neighborhood $ \mcV= [t_1,t_2]\times V $ of a given point~$m_0$, and
choose coordinates $ x^j$ on~$V $; this provides coordinates $(t,x^i)$
on $\mcV$, which induce coordinates on $T^1\mcV$: for $ m\in\mcV$
and $ {\dot m}\in T^1_m\MMM$, write  ${\dot m} = {\dot m}^0\,\partial
_t + \sum_{1\le j\le d} {\dot m}^j\,\partial_{x^j}$.

Note first that since $ U = a^{-2} \partial_t $, we have
\[
\nabla_{\dot m}U = \nabla_{\dot m}(a^{-2})\, \partial_t+a^{-2}\nabla
_{\dot m}\,\partial_t .
\]
Using Christoffel's symbols $\Gamma^i_{jk}$ we have
\[
(\nabla_{\dot m}\,\partial_t )^{\al} = \nabla_{\dot m}(a^{-2}) \delta
_0^{\al} + a^{-2} {\dot m}^c \Gamma^{\al}_{c 0},
\]
for $\al\in\{0,\ldots,d\}$ and a summation over $c$ in $\{0,\ldots,d\}$; so
\[
\label{EqIntermediaire}
H_0 f = g(\nabla_{\dot m}U,\dot m) = \nabla_{\dot m}(\log a^{-2})
{\dot m}^0 + a^{-2}{\dot m}^c \Gamma^{\al}_{c 0} g_{\al\beta} {\dot
m}^{\beta}.
\]
The explicit formulas for the Christoffel symbols, in terms of the
metric, are
\begin{eqnarray*}
\Gamma^0_{0 0} &=& \partial_t(\log a),\qquad
\Gamma^0_{k 0} = \partial_{x^k}(\log a),\\
\Gamma^i_{0 0} &=& \tfrac{1}{2} h^{i \ell}\,\partial_{x^{\ell}}(a^2),\qquad
\Gamma^i_{k 0} = \tfrac{1}{2} h^{i \ell}\,\partial_th_{\ell k} ,
\end{eqnarray*}
for $ i,k\in\{1,\ldots,d\}$ and a sommation over $ 1\le\ell\le d $. We
thus have, after simplifications,
\begin{eqnarray*}
H_0f &=& -2 \nabla_{\dot m}(\log a) {\dot m}^0 + |{\dot m}^0|^2
\,\partial_t(\log a) - \frac{a^{-2}}{2} {\dot m}^k \,\partial_t(h_{\ell
k}) {\dot m}^{\ell} \\
&=& - |{\dot m}^0|^2 \,\partial_t\log a - 2 {\dot m}^0 {\dot m}^k
\,\partial_{x^k}\log a - \frac{a^{-2}}{2} {\dot m}^k \,\partial
_t(h_{\ell k}) {\dot m}^{\ell}.
\end{eqnarray*}
Using the unit pseudo-norm relation $ a_m^2 |\dot m^0|^2-h_{\ell k}(m)
{\dot m}^k {\dot m}^{\ell}=1$, the above equality becomes$ $
\[
H_0f = -|\dot m^0|^2 \,\partial_t\log a - 2 {\dot m}^0{\dot m}^k
\,\partial_{x^k}\log a - \frac{a^{-2}}{2} |\dot m^0|^2 \,\partial_t(a^2) ,
\]
that is,  $H_0f = -2 {\dot m}^0 \nabla_{\dot m}\log a $.  The
statement of the theorem follows from~\eqref{CriterionNonExplosion}.
\end{pf}

This result takes a particularly simple form in the case where $\Theta
$ depends only on the base point $ m $, as is the case of the $R$-diffusion.\vadjust{\goodbreak}

\begin{cor} Let $ \MMM= I\times S $ denote a generalized warped product spacetime
and $ \Theta$ be a bounded nonnegative function on $\MMM$. Then the
$\Theta$-diffusion does not explode if $ \nabla a $ is everywhere
nonspacelike and future-directed.
\end{cor}

\begin{pf}
$\!\!\!$The condition of Theorem~\ref{ThmWarped} reads, in that case,
``\textit{$T^1\MMM\ni(m,\dot m)\mapsto\nabla_{\dot m} \log a $ is
bounded below}.'' To rephrase this condition into the more synthetic
condition of the statement, let us work in local coordinates, $(t,x)$
and $(\dot t, \dot x)$ for $m$ and $\dot m$, respectively.

We have $ \dot t = a^{-1}\operatorname{ch}r $ and $ \dot x = (\operatorname{sh}r)\sigma$, for some $r\in\RR$ and $\sigma\in T_xS$ with
$|\sigma|_{h( m)} = 1$.

Define $ u:= \partial_t\log a $ and $ v := \partial_x\log a\in
T_xS\equiv\RR^d$. Then the condition of Theorem~\ref{ThmWarped}
reads, $ $ ``$u a^{-1} \operatorname{ch} r - (v_i \sigma^i) \operatorname{sh} r
\geq C$,'' for any $r$ and $\sigma$. Letting $r\ra\pm\infty$, gives
$ a^{-1} u\ge|v_i \sigma^i| \ge0 $. As the constant $C$ can be taken
negative without loss of generality, the reciprocal is clear. Now,
since $\max_{|\sigma|_{h( m)}=1} |v_i \sigma^i| = |v|_{h^{-1}( m)}
$, the condition reads$ $,
``$a^{-1} u\ge|v|_{h^{-1}(m)} $.'' Finally, as $\nabla= (
a^{-2}\,\partial_t ,\break -h^{ij}\,\partial_{x^j} )$, the vector $\nabla\log a
= ( a^{-2}u ,-h^{ij}v_j )$ has pseudo-norm $ g (\nabla\log a,\break\nabla
\log a ) = a^{-2} u^2 - |v|^2_{h^{-1}( m)}\ge0 $.
\end{pf}

 This criterion applies in particular to $\Theta$-diffusions in
Robertson--Walker spacetimes, recovering the results of Angst~\cite{JurgenPhD}, who proceeded by direct analysis of the stochastic
differential equations of the dynamics.

%%----------------------%%
%s3.2 ###
\subsection{Perfect fluids}
\label{SectionPerfectFluids}
%%----------------------%%

Our second class of examples, where to apply Lyapounov's method to
prove nonexplosion, will be the set of spacetimes with normal matter
whose energy-momentum tensor $\bfT$ is that of a perfect fluid. They
are characterized by the datum of a timelike vector field $ U$, the
four velocity of the fluid and two functions $\rho$ and $p$ on $\MMM
$, respectively, the energy density and pressure of the fluid. See \cite
{HawkingEllis,BeemEhrlich}. We have then $\bfT= \rho
U\otimes U+ p (g + U\otimes U )$, or in local coordinates,
\[
\bfT_{ij} = (\rho+p) U_iU_j + p g_{ij} .
\]
Such a spacetime is said to be of \textit{perfect fluid
type}. Notice that contrarily to the globally hyperbolic spacetimes,
no topological assumption is made on a~perfect fluid type spacetime.

G\"{o}del's universe is such a spacetime. This is the manifold $\RR^4$
with the metric $ds^2 = dt^2 -dx^2+\frac{1}{2}e^{2\sqrt{2}\omega x}
dy^2-dz^2-2e^{\sqrt{2}\omega x} \,dt \,dy$, where $\omega>0$ is a~constant.
It is a solution to Einstein's equation with cosmological
constant $\omega^2$ and represents a pressure-free perfect fluid. It
has energy-momentum tensor  $\mathbf{T} = U\otimes U $,
where $(U_j)=(\sqrt{2}\omega,0,\sqrt{2}\omega e^{\sqrt{2}\omega
x},0)$ represents the four-velocity covector of the matter, and $\omega
$ is the vorticity of this field. This spacetime has constant scalar
curvature $ 2\omega^2$. See Section 2.4 in~\cite{FranchiGodel}. As
above, the function~$f$ is defined by formula~(\ref{Definitionf})
and can be used as a Lyapounov function under some conditions. The
computations made in Section~\ref{SectionGloballyHyperbolic} work
equally well in that setting and lead to the following results.\vadjust{\goodbreak}

\begin{prop}
Let $(\MMM,g)$ be a Lorentzian manifold of perfect fluid type and $ f
$ be defined by formula~(\ref{Definitionf}). Suppose $f$ goes almost-surely to
infinity along any exploding timelike path. If there exists a constant
$ C $ \st
\[
H_0f + \frac{d}{2} \Theta f+ \frac{1}{2} (f \dot m^k-U^k
)\frac{\partial\Theta} {\partial\dot m^k} \leq C f ,
\]
then the $\Theta$-diffusion has almost-surely an infinite lifetime.
\end{prop}

In the particular case of G\"{o}del's universe, the gradient $\nabla U$
of the velocity vanishes (since $U^i = \delta^{i}_0$), so that $H_0f =
0$, by formula~\eqref{FormulaH_0f}; and $f$ is the square root of the energy.

\begin{cor} \label{cor.Godel}
Let us work in G\"{o}del's universe and suppose that $ 3\Theta+ (\dot
m^k\frac{\partial\Theta} {\partial\dot m^k} - \frac{1}{f}\frac
{\partial\Theta}{\partial{\dot m}^0} )$ is bounded above in $T^1\MMM
$. Then the $ \Theta$-diffusion has almost-surely an infinite lifetime.
This condition holds in particular if $ \Theta(\Phi)=\Theta(m)$
depends only on the base point and is bounded, as this is the case for
the basic relativistic diffusion and the $ R$-diffusion in G\"{o}del's universe.
\end{cor}

Note that this criterion does not apply to the energy diffusion in G\"
{o}del's universe. Indeed one can see in that case (see Section 2.4 of
\cite{FranchiGodel}) that the above quantity is equal to $ 5\Theta
-4\omega^2$ and that the energy $ \Theta$ is unbounded along the
trajectories of the energy diffusion.

\begin{rem}
In Einstein--de Sitter spacetime the energy diffusion explodes
with positive probability, as proved in Proposition 5.4.2 of \cite
{FranchiLeJanCurvature}. (This Robertson--Walker universe is both a
warped product and a perfect fluid type spacetime.) Consult \cite
{BailleulJMP} for a first study of stochastic incompleteness for
relativistic diffusions.
\end{rem}

%--------------------------%
%s4 ###
\section{b-completeness}
\label{SectionbCompleteness}
%--------------------------%

The study of dynamics in the orthonormal frame bundle is not new in
general relativity, and essentially dates back to Cartan's moving frame
method. However, Schmidt~\cite{SchmidtbComplete1} was the first to
notice that the geometry of $\OO\MMM$ itself may be used to provide a
conceptual framework in which studying the nature of spacetime
singularities. For that purpose, he introduced on the parallelizable
manifold $\OO\MMM$ a Riemannian metric, turning $\{H_0,\ldots,H_d
,(V_{ij})_{0\leq i<j\leq d}\}$ into a Riemannian orthonormal basis, and
called it the bundle metric, or \textit{$ b$-metric}.
The completeness of this metric structure on $\OO\MMM$ can
essentially be phrased in terms of $\MMM$-valued paths. To state that
fact, recall that one can associate to any $\MMM$-valued $\mcC^1$
path $\gamma\dvtx [0,T[\, \ra\MMM$ and $\be\in\OO_{\gamma_0}\MMM$ a
unique horizontal lift $ \gamma^{\uparrow} \dvtx  [0,T) \ra\OO\MMM$ of
$\gamma$, starting from $(\gamma_0,\be)$, and charactarized by the properties
\[
\frac{d}{ds}\gamma^{\uparrow}_s\in\operatorname{span}(H_0,\ldots,H_d)
\quad \mbox{and} \quad \pi_0 (\gamma^{\uparrow}_s )=\gamma_s \qquad \mbox{ for
all } s\in[0,T) .
\]
The \textit{$S_\be$-length} of $ \gamma$ is defined as the
Riemannian length of its horizontal lift $ \gamma^{\uparrow}$; it
depends on\vadjust{\goodbreak} $\be\in\OO_{\gamma_0}\MMM$. In other words, given $\be
\in\OO_{\gamma_0}\MMM$, seen as orthonormal in the Euclidean sense,
the $S_\be$-length of the $\MMM$-valued $\mcC^1$ path~$\gamma$ is
the Euclidean length of its anti-development in $ (T_{\gamma_0}\MMM
,\be)$. Although this length depends on $\be$, its finiteness is
independent of it; we can thus talk of \textit{finite $S$-length} of a
$\mcC^1$ path without mentioning the frame $\be$. Note that in
a~Riemannian setting the $S_\be$-lenth of a $\mcC^1$ path is its usual
Riemannian length.

\begin{thm}[(Schmidt~\cite{SchmidtbComplete1})]
$\OO\MMM$ is complete for the above b-metric if and only if  any $\mcC^1$ path
$\gamma\dvtx [0,T)\ra\MMM$ with a bounded $S$-length converges in~$\MMM$ at time $T^-$.
\end{thm}

 The above completeness hypothesis is usually called \textit{b-completeness}. The Riemannian version of this statement is
trivial as the orthonormal frame bundle with its b-metric is complete
iff the Riemannian manifold is complete. The Lorentzian situation is
more involved as there exists (timelike, spacelike and lightlike)
complete Lorentzian manifolds $\MMM$ which have an incomplete path of
bounded acceleration, so $\OO\MMM$ is not b-complete (see, e.g.,
\cite{GerochSingularity} and~\cite{BeemCounterexample}). The
noncompactness of $\mathit{SO}_0(1,d)$ lies at the core of this phenomenon.

However, the Riemannian view of a Lorentzian manifold provided by
Schmidt's metric offers a bridge to investigate some features of the
latter using the tools of Riemannian geometry, as the following
proposition shows.

\begin{prop}\label{prop}
Let $\Theta$ be a bounded function on $\MMM$. Then the $\Theta
$-diffusion does not explode if $\OO\MMM$ is b-complete.
\end{prop}

One should not be confused about that statement. It does not mean that
the Riemannian completeness of $\OO\MMM$ implies the completeness of
its Brownian trajectories, which is false. One cannot assign an $S_\be
$-length to a Brownian path in $\OO\MMM$ as it is not regular enough.

\begin{pf*}{Proof of Proposition~\ref{prop}}
$\bullet$ Given a horizontal $\mcC^1$-path $(\rho_s)_{0\leq s<T}$ in
$\OO\MMM$, write $\gamma$ for its projection $\pi_0\circ\rho$ in
$\MMM$, so $\rho=\gamma^{\uparrow}$. For $0\le s<T$, denote by $
\tau^\gamma_{0\to s} $ the parallel transport operator along the
curve $(\gamma_r)_{0\leq r\leq s}$, with inverse $ \tau^\gamma
_{0\lar s} $. Also, denote by $(p_s)_{0\leq s<T} $ the anti-development
of $ \gamma$: this $ T_{\gamma_0}\MMM$-valued $\mcC^1$-path is
defined for all $ s\in[0,T[ $ by the formula$ $ $ {  p_s = \int_0^s
\tau^\gamma_{0\leftarrow r} \dot\gamma_r \,dr} $. Last, we shall
denote by $ {\dot{p}}^j_r $ the coordinates of $ \dot{p}_r $ in the
frame $\rho_0$, and by $\| . \|_{\rho_s} $ the Euclidean norm in
$(T_{\gamma_s}\MMM,\rho_s)$.  We have by construction $d\rho_s
= \sum_{0\le j\le d} H_j(\rho_s) {\dot{p}}^j_s \,ds $ and $ \dot
{\gamma}_s = \tau_{0\to s}^\gamma\dot{p}_s $, as well as the
identity $\|\dot\gamma_s\|_{\rho_s}^2 = \|\dot p_s\|_{\rho_0}^2 =
\sum_{0\le j\le d} ({\dot{p}}^j_s )^2$. The b-completeness assumption\vspace*{2pt}
means that $\gamma$ has a limit $\gamma_T$ in $\MMM$ at time $T$ if
%
%e4.1 ###
\begin{equation}
\label{ConditionSchmidt}
\int_0^{T-} \|\dot p_s\|_{\rho_0} \,ds < \infty.\vadjust{\goodbreak}
\end{equation}

$\bullet$ The basic relativistic diffusion $ (m_s,\be_s )_{0\leq
s<\zeta} $ is by construction the development in $\MMM$ of the
relativistic Dudley diffusion in Minkowski spacetime, identified with
$T_{m_0}\MMM$ (see Theorem 3.2 in~\cite{FranchiLeJan}). As trajectories of the
latter over a time a bounded time interval have almost-surely a finite length in
the Eulidean norm associated with any frame of $\RR^{1,d}$, the
b-completeness of~$\OO\MMM$ ensures the nonexplosion of the basic
relativistic diffusion.

$\bullet$ For a generic $\Theta$-diffusion, formula \eqref
{GeneratorCoordinates} implies the existence for each $s\in[0,\zeta[$
of an orthonormal basis $ (\varphi_1(s),\ldots,\varphi_d(s) )$ of $ \dot
{p}_s^\bot$ in $\RR^{1,d}$ \st one has
\[
d\dot p_s^k = \sum_{j=1}^d \sqrt{\Theta(m_s)} \varphi^k_j(s)
\,dw^{j}_s + \frac{d}{2} \Theta(m_s) \dot p_s^k \,ds
\]
for some $d$-dimensional Brownian motion $w$. We have used the fact that $\Theta$
depends only on $m$ to simplify the general expression. The path
$(p_s,\dot p_s)_{0\leq s<\zeta}$ appears then as a time change of
Dudley's diffusion, by means of the map $s\mapsto\inf\{u | \int
_0^u\Theta(m_r) \,dr>s \}$. The result follows for a bounded function
$\Theta$.~%
\end{pf*}

This result can be improved in two ways: by relaxing the boundedness
hypothesis on $\Theta$ and by relaxing the geometric completeness
assumption. The next section explains how this can be done in a
sub-Riemannian framework by using ideas from the theory of reversible
Markov processes.

%-----------------------------------------------------%
%s5 ###
\section{A volume growth nonexplosion criterion}
\label{SectionVolumeGrowth}
%-----------------------------------------------------%

We prove in this section a nonexplosion criterion involving only the
volume growth of some sub-Riemannian boxes in $\OO\MMM$ and the
function $\Theta$, as described in Theorem~\ref{MainThm} below. This
result is proved in Section~\ref{SectionCrossingTimes} following
Takeda's method, as improved recently by Hsu and Qin in~\cite{HsuQin}.
Yet, there is a real difficulty in doing this, as we are working with a
nonsymmetric, hypoelliptic diffusion, and on a principal bundle with
noncompact fibers. To overcome these difficulties, we introduce a
sub-Riemannian structure on $\OO\MMM$, well adapted to our setting,
and which will somehow play for us the role of the missing Lorentzian distance.

%%---------------------------------------------------%%
%s5.1 ###
\subsection{Sub-Riemannian framework and main results}
%%---------------------------------------------------%%

%%%--------------------------------------%%%
%s5.1.1 ###
\subsubsection{Sub-Riemannian distance function}
\label{SubsectionSubRiemannian}
%%%--------------------------------------%%%

We have seen in Section~\ref{SectionbCompleteness} that the
completeness of the natural Riemannian metric of the parallelizable
manifold $\OO\MMM$ implies the stochastic completeness of all the
$\Theta$-diffusions with a bounded $\Theta$. One can significantly
improve that conclusion by working with the sub-Riemannian structure on
$\OO\MMM$ induced by the field of $(d+1)$-planes generated by the
vector fields $H_0,V_1,\dots,V_d$. In that setting, one can assign a
length only to $\mcC^1$ paths $\rho\dvtx  [0,T]\ra\OO\MMM$ whose
tangent vector belong at any time $s$ to the vector space spanned by
$H_0,V_1,\dots,V_d$ in $T_{\rho_s}\OO\MMM$, say $\dot\rho_s =
\dot\rho_s^0 H_0(\rho_s) + \dot\rho_s^1 V_1(\rho_s) + \cdots+
\dot\rho_s^d V_d(\rho_s)$. Such a path is said\vadjust{\goodbreak} to be \textit
{admissible}; its length is then defined as $\int_0^T (\sum_{i=0}^d
(\dot\rho_s^i )^2 )^{{1}/{2}} \,ds$. The sub-Riemannian distance
between two points of $\OO\MMM$ is defined as the infimum of the
length of the admissible paths joining these two points, with the
convention $\inf\emptyset= +\infty$. Chow's theorem~\cite{Chow}
ensures that the sub-Riemannian distance function $\mcD(\cdot,\cdot
)$ is finite and continuous in its two arguments if (see, e.g., \cite
{Montgomery}) the Lie algebra generated by $H_0,V_1,\dots,V_d$ has
full dimension, which holds here. Fix a reference point $\Phi_{\mathrm{ref}}\in\OO\MMM$.

\begin{HCH*}\hypertarget{hypH}{}
$\!\!\!$The closed boxes
\mbox{$B_\la:=\{\mcD(\Phi_{\mathrm{ref}},\cdot)\leq\la\}$} are compact
for any $\la>0$.
\end{HCH*}

This completeness hypothesis rules out the pathological examples of
Geroch~\cite{GerochSingularity} and Beem~\cite{BeemCounterexample};
it does not depend on the arbitrary choice of $\Phi_{\mathrm{ref}}$.
Unlike its Riemannian analog, the sub-Riemannian distance function
$\mcD(\Phi_{\mathrm{ref}},\cdot)$ is not smooth in any neighborhood
of $\Phi_{\mathrm{ref}}$,~\cite{Montgomery}; however, it is a
viscosity solution of the equation
\[
|H_0\mcD|^2 + |V_1\mcD|^2 + \cdots+ |V_d\mcD|^2 = 1
\]
on $\OO\MMM\bs\{\Phi_{\mathrm{ref}}\}$ (see, e.g., Theorem 2 in
\cite{Dragoni07}; we do not use that fact in the sequel). We shall use
that quantitative information in Section~\ref{SectionCrossingTimes}
under the classical form given in the following proposition.

\begin{prop}
\label{CorApproxTimeFunction}
Fix $\la>0$. One can associate to any positive constant~$\eta$ a
smooth function $F \dvtx  \OO\MMM\ra\RR_+$ \st
\[
\max_{\Phi\in B_\la} |F(\Phi)-\mcD(\Phi_{\mathrm{ref}},\Phi) | \leq\eta
\]
and we have on $B_\la$
\[
|H_0F|^2 + |V_1F|^2 + \cdots+ |V_dF|^2 \leq2.
\]
\end{prop}

\begin{pf}
Let us introduce the Riemannian metric $g_\ep$ on $\OO\MMM$ for
which $H_0,H_1,\dots, H_d$ and the $(V_{ij})_{0\leq i<j\leq d}$ are
\textit{orthogonal}, with $H_0$ and the \mbox{$V_{0j}(=V_j)$} of norm $1$ and
the other vectors of norm $\ep^{-1}$. Denote by $\mcD_\ep(\cdot) =
\mcD_\ep(\Phi_{\mathrm{ref}},\cdot)$ the distance function
associated with $g_\ep$. It is a $1$-Lipschitz-continuous function
(\wrt the distance function $\mcD_\ep$) which is differentiable
almost-everywhere, by Rademacher's theorem, and has a gradient of norm
$1$ almost-everywhere.
%
%e5.1 ###
\begin{eqnarray}
\label{RademacherEq}
&&|H_0\mcD_\ep|^2 + |V_1\mcD_\ep|^2 + \cdots+ |V_d\mcD_\ep|^2 \nonumber\\[-8pt]\\[-8pt]
&&\qquad {}+
\ep^{-2} \Biggl(\sum_{i=1}^d|H_i\mcD_\ep|^2 + \sum_{1\leq i<j\leq d}
|V_{ij}\mcD_\ep|^2 \Biggr)= 1.\nonumber
\end{eqnarray}
(Indeed, the set of conjugate points to $\Phi_0$ in $B_\la$ is closed
and has null measure. In the complementary, relatively open, set the
distance is attained along a unique geodesic whose unit tangent vector\vadjust{\goodbreak}
at the final point is the gradient of the distance function to $\Phi
_0$.) The function $\mcD_\ep$ is easily seen to converge uniformly to
$\mcD(\Phi_{\mathrm{ref}},\cdot)$ on the compact box $B_\la$ (this
is where we need these boxes to be compact); see, for example, Sections
0.8.A and 1.4.D of Gromov's article~\cite{GromovCC}. As we have
almost-everywhere
\[
|H_0\mcD_\ep|^2 + |V_1\mcD_\ep|^2 + \cdots+ |V_d\mcD_\ep|^2 \leq1,
\]
by~\eqref{RademacherEq}, a standard regularization procedure yields
the conclusion.
\end{pf}

%%%--------------------%%%
%s5.1.2 ###
\subsubsection{Main results}
%%%--------------------%%%

We use the natural volume measure on $\OO\MMM$ associated with the
Lorentzian structure. It is defined by the formula
\[
\Vol(d\Phi) = \Vol_\MMM(dm)\otimes\Vol_m(d\be),\qquad  \Phi=(m,\be),
\]
where $\Vol_\MMM(dm)$ is the Lorentzian volume measure and $\Vol
_m(d\be)$ is the image of a given Haar measure on $\mathit{SO}_0(1,d)$ by the
identification of the fiber $\pi_0^{-1}(m)$ with $\mathit{SO}_0(1,d)$; see, for
example,~\cite{HawkingEllis}, Section 2.8, for the Lorentzian volume
measure. The volume measure $\Vol$ on $\OO\MMM$ is uniquely defined
up to a multiplicative constant. In order to avoid some unpleasant
pathologies, we shall make the following rather mild assumption on the
causal structure of spacetime.

\begin{hyp*}
$(\MMM,g)$ is strongly causal.
\end{hyp*}

 It means that any point of $\MMM$ has arbitrarily small
neighborhoods which no nonspacelike path intersects more than once; see
\cite{HawkingEllis}, page 192, or~\cite{BeemEhrlich}.

\begin{thm}
\label{MainThm}
Let $(\MMM, g)$ be a strongly causal Lorentzian manifold satisfying
the Completeness hypothesis \hyperlink{hypH}{(\textup{H})}. Set $\Theta_r :=
\sup_{\Phi\in B_r } \Theta(\Phi)$, for any $r>0$, and suppose
%
%e5.2 ###
\begin{equation}
\label{ConditionMainThm}
\int^\infty\frac{r \,dr}{\Theta_r \log(\Theta_r \Vol(B_r) )} =
\infty .
\end{equation}
Then the $\Theta$-diffusion has almost-surely an infinite lifetime, from any
starting point.
\end{thm}

Condition~\eqref{ConditionMainThm} has the form of the classical
nonexplosion condition for Brownian motion, $  \int^\infty \frac
{r \,dr}{\log\Vol(B_r)} = \infty$, first proved by Grigor'yan \cite
{GrigoryanCompleteness} and has precisely that form for $\Theta$
bounded. Note that no topological assumption on $\MMM$ is needed,
contrary to the results of Section~\ref{SectionGloballyHyperbolic}.
One can give a quantitative version of the above theorem by providing
an upper rate function.

\begin{cor}
\label{CorMainThm}
Let $ \MMM$ be a strongly causal Lorentzian manifold satisfying the
Completeness hypothesis \hyperlink{hypH}{(\textup{H})}. Set $h(\rho) \equiv
\rho$ if $\Theta\equiv0$; otherwise, pick a constant $R_0$ such that
$\Theta_{R_0}>0$ and set for $\rho>0$
\[
h(\rho) := \inf\biggl\{R>R_0 \Big| \int_{R_0}^R \frac{r \,dr}{\Theta_r \log
[\Theta_r \Vol(B_r) ] } > \rho\biggr\}.\vadjust{\goodbreak}
\]
Then, given any $ \Phi_0\in\OO\MMM$, there exist $R_0>0$ and a
positive constant $C$ \st we have $\PP_{\Phi_0}$-almost-surely
\[
\mcD(\Phi_0,\Phi_s) \leq C h(Cs).
\]
\end{cor}

We prove Theorem~\ref{MainThm} following Takeda's method, explained in
the next section. To adapt it to our setting, we shall introduce in
Section~\ref{SectionModifiedProcess} a modified $\Theta$-diffusion on
some compact space; it is used crucially in the proof of Theorem \ref
{MainThm} given in Section~\ref{SectionCrossingTimes}.

%%-------------------------%%
%s5.2 ###
\subsection{Takeda's method}
\label{SectionMethod}
%%-------------------------%%

%%%---------------------------%%%
%s5.2.1 ###
\subsubsection{The main ingredients}
%%%---------------------------%%%

Using an idea of Lyons and Zheng~\cite{LyonsZheng}, Takeda devised
\cite{Takeda1,Takeda2} a remarkably simple and sharp
nonexplosion criterion for Brownian motion  on a Riemannian manifold $\VV$. Loosely
speaking, his reasoning works as follows. Suppose we have a diffusion
$(x_s)_{s\geq0}$ on $\VV$ which is symmetric (\wrt the Riemannian
volume measure $\Vol$, say) and conservative; denote by $L$ its
generator, and let $f$ be a sufficiently smooth function. Denote by $
\PP_{\Vvol}$ the measure $\int\PP_x \Vol(dx)$ on the path space,
where $ \PP_x$ is the law of the diffusion started from $x$. Fix a
time $T>0$. As the reversed process $(x_{T-s})_{0\leq s\leq T}$ is an
$L$-diffusion under $\PP_{\Vvol} $, applying It\^{o}'s formula to both\vspace*{2pt}
$f(x_s)$ and $f(x_{T-s})$ provides two martingales $M$ and $\widetilde
{M}$ [\wrt the two different filtrations $\sigma(x_s ; 0\leq s\leq T)$
and $\sigma(x_{T-s} ; 0\leq s\leq T)$, resp.] \st
\begin{eqnarray*}
f(x_s) &=& f(x_0) + M_s + \int_0^sLf(x_r) \,dr , \\
f(x_s) &=& f\bigl(x_{T-(T-s)}\bigr) = f(x_T) + \widetilde{M}_{T-s} + \int_0^{T-s}Lf(x_{T-r}) \,dr.
\end{eqnarray*}
It follows that $f(x_s) = \frac{f(x_0)+f(x_T)}{2} + \frac
{M_s+\widetilde{M}_{T-s}}{2} +\int_0^TLf(x_s) \,dr$, and consequently,
\[
f(x_T)-f(x_0) = \tfrac{1}{2} (M_T-\widetilde{M}_T ).
\]
If $\frac{d\langle M\rangle_s}{ds}$ and $\frac{d\langle\widetilde
{M}\rangle_s}{ds}$ are bounded above, by $1$ say, the previous
identity provides a control of $ (f(x_T)-f(x_0) )$ by the supremum of
the absolute value of a Brownian motion  over the time interval $[0,T]$.

Back to the nonexplosion problem for Brownian motion  on $\VV$, fix a point \mbox{$ m\in
\VV$} and a radius $R>1$, and consider the Brownian motion $(x_s)_{s\geq0}$
reflected on the boundary of the Riemannian ball $B(m;R)$, started
under its invariant measure $\mathbf{1}_{B(m;R)}\Vol$. It is a symmetric
conservative diffusion$ $; denote by $\overline\PP_{B(m;R)}$ its law.
Using the Dirichlet forms approach to symmetric diffusions one can
apply the above reasoning to the (nonsmooth,\vadjust{\goodbreak} but $1$-Lipschitz)
Riemannian distance function $d(m,\cdot)$, which gives the estimate
\[
\overline\PP_{B(m;R)} \Bigl(x_0\!\in\! B(m;1), \sup_{s\leq T}d(m,x_s)\!=\!R \Bigr)
\!\leq\!\Vol(B(m;R) )\!\times\!2 \PP\Bigl(\sup_{s\leq T}|B_s|\!>\!R \Bigr).
\]
But as the Brownian motion  on $\VV$ behaves in the ball $B(m;R)$ as the Brownian motion
reflected on the boundary of $B(m;R)$, the above inequality also gives
an upper bound for the probability that the Brownian motion  on $\VV$, started
uniformly from $B(m;1)$, exits the ball $B(m;R)$ before time $T$.
Combining this estimate with the Borel--Cantelli lemma, Takeda proved
that the Brownian motion  on $\VV$ is conservative provided
\[
\liminf_{R\ra\infty} {R^{-2}} \log\Vol(B(m;R) )<\infty,
\]
re-proving in a simple way a criterion due to Karp and Li. Takeda's
method has been refined by several authors, culminating with Hsu and
Qin's recent work~\cite{HsuQin}, in which they give an elegant and
simple proof of a sharp nonexplosion criterion, due to Grigor'yan \cite
{GrigoryanCompleteness}, for Brownian motion  on a~Riemannian manifold in terms of
volume growth, as well as an escape rate function. We shall follow
their method to deal with relativistic diffusions.

%%%----------------------%%%
%s5.2.2 ###
\subsubsection{The difficulties}
\label{SubSectionDifficulties}
%%%----------------------%%%

The main difficulty in implementing this approach is in finding what
can play the role of the pair ``\textit{Riemannian distance function--reflected Brownian motion}'' in our Lorentzian, hypoelliptic
framework. We describe in the remainder of this section a nonstandard
reflection mechanism for a Brownian motion  in a Riemannian manifold which will
serve us as a~guide in the construction of the $\Theta$-diffusion
reflected on the boundary of the sub-Riemannian boxes, as described in
Section~\ref{SectionModifiedProcess}.

Brownian motion  reflected on the boundary of a ball $B(m;R)$ is the simplest
diffusion process which coincides with Brownian motion  on the ball $B(m;R)$ and has
a state space with finite volume. One cannot take a~smaller state space
if the former property is to be satisfied. Yet, one can make different
choices if one is ready to loose the minimality property. To explain
that fact, let us suppose that $(\VV,g)$ is a Cartan--Hadamard
manifold. Given a point $m\in\VV$ let us use the exponential map
$\exp_m$ at $m$ as a global chart on $\VV$; this identifies the
geodesic ball $B(m;R)$ on $M$ to the (Euclidean-shaped) ball $B'(0;R)$
in $T_m\VV$. Given $\ep>0$, let us modify the metric on $B'(0;R+\ep
)\moins B'(0;R)$ so as to interpolate smoothly between $\exp_m^*g$ on
$B'(0;R)$ and the constant metric $g_m$ outside $B'(0;R+\ep)$ (primed
balls refer to the pull-back metric $\exp_m^*g$). Denote by
$\widetilde{g}$ the restriction to $B'(0;R+2\ep)$ of this modified
metric, and define the compact space $\mathbf{K}$ as the quotient of the
closed ball $\overline{B}'(0;R+2\ep)$ by the identification of $m'\in
\partial\overline{B}'(0;R+2\ep)$ and $-m'$. Then the $\widetilde
{g}$-Brownian motion  on $\mathbf{K}$ coincides with the $\exp_m^*g$-Brownian motion  on $B'(0;R)$
and has a state space with finite $\widetilde{g}$-volume $\Vol
_{\widetilde{g}}(\mathbf{K}) = (1+o(\ep)) \Vol_g (B(m;R) )$. The
construction of a modified\vadjust{\goodbreak} $\Theta$-diffusion given in Section \ref
{SectionModifiedProcess} will be reminiscent of the preceding
nonstandard reflected Brownian motion.

%%%-------------------------%%%
%s5.3 ###
\subsection{A modified process}
\label{SectionModifiedProcess}
%%%-------------------------%%%

We start our construction of the ``reflected'' $\Theta$-diffusion by
constructing the compact space on which it is going to live. Fix for
that purpose a reference point $\Phi_{\mathrm{ref}} \in\OO\MMM$,
the center of the boxes $B_\la$, and set $\mcD(\Phi) = \mcD(\Phi
_{\mathrm{ref}},\Phi)$ for all $\Phi\in\OO\MMM$. Fix also two
positive constants~$\la$ and $\varep$ and consider the relatively
compact open region
\[
\mcU:= \{\la< \mcD<\la+\varep\} = B_{\la+\ep} \bs\overline
{B_\la}.
\]

\begin{lem}
\label{lem.hyps}
There exists in $\mcU$ a smooth hypersurface $V$ of $\OO\MMM$
separating $\partial B_{\la}$ from $\partial B_{\la+\varep}$ \st the
subset $V_0 := \{\Phi\in V | H_0(\Phi)\in T_\Phi V\}$ is a~smooth
hypersurface of $V$.
\end{lem}

 The separation property means that $ \partial B_{\la}\cup
\partial B_{\la+\varep}$ does not intersect $V$ but any continuous
path from $\partial B_{\la}$ to $\partial B_{\la+\varep}$ hits $V$.
We thank A. Oancea and P. Pansu for their help in proving this statement.

\begin{pf*}{Proof of Lemma~\ref{lem.hyps}}
Let us use the function $F$ of Proposition \ref
{CorApproxTimeFunction}, with $ \eta<\varep/4 $ and $ R>\la+\varep
$, and fix some constants $ \eta< \varep_1<\varep_2<\varep/2 -\eta
$ such that $B_{\la}\subset\{\varep_1\le F-\la\le\varep_2\}
\subset B_{\la+\varep/2} $. The set of regular values of $(F-\la)$
is dense in the interval $(\varep_1,\varep_2) $, by Sard's theorem.
Fix a regular value $ c\in(\ep_1,\ep_2)$, so the level set $ S:= \{
F= c\}$ is a smooth hypersurface separating $ \partial B_{\la} $ from
$ \partial B_{\la+\varep/2} $.

We shall now be working in $ \mcU' \equiv S \times [0,\frac{\varep
}{2} )$, where we are going to construct the separating hypersurface
$V$ as the graph of some function $f \dvtx  S \ra [0,\frac{\varep}{2} )$,
resorting to the transversality lemma. Denote by $\operatorname{Gr}(T\mcU')$ the
Grassmannian bundle over $\mcU'$ made up of all the hyperplanes of
$T\mcU'$, and associate to any function $ f \dvtx  S\to (0,\frac{\varep
}{2} )$ the function $ G_f \dvtx  S\to \operatorname{Gr}(T\mcU')$ defined by $ G_f(m) :=
\{ (\sigma, df_m(\sigma) ) | \sigma\in T_mS \}$. Let $\mcH$ denote
the smooth hypersurface of $\operatorname{Gr}(T\mcU')$, made up of all hyperplanes
containing $H_0 $. Then $ G_f^{-1}(\mcH)$ is a~smooth hypersurface of
$ \operatorname{Graph}(f)$ as soon as $G_f$ is transverse to $\mcH$. Therefore
the statement reduces to finding a function $ f $ such that $G_f$ be
transverse to $\mcH$.

Consider for that purpose a smooth partition of unity: $\mathbf{1}_S=\sum
_{j=1}^k \alpha_j $, with $\{\alpha_j>0\}=\psi_j(\mcB^\nu)$
diffeomorphic under $\psi_j$ to the unit ball $\mcB^\nu\subset\RR
^\nu$ [with $ \nu=\operatorname{dim}(\OO\MMM)-1=(d+3)d/2$]. Denoting by
$\mcA$ the space of (the restictions to $\mcB^\nu$ of) affine
functions on $\RR^\nu$, consider the map $ F \dvtx  \mcA^n\times S\to
\operatorname{Gr}(T\mcU')$ defined by the formula
\[
G(\varphi_1,\ldots, \varphi_k,m) := G_f(m),
\]
where $f = \sum_{j=1}^k \alpha_j \varphi_j\circ\psi_j^{-1}$. This
is easily seen to be a submersion. It follows from the transversality
lemma that such a $G_f$ is transversal\vadjust{\goodbreak} to $ \mcH$ for almost-every
$(\varphi_1,\ldots, \varphi_k)\in\mcA^n$. The graph of the
function $f$ corresponding to a~small multiple of such a $k$-tuple has
the properties of the statement.~%
\end{pf*}

Let $O$ be the set of points of the box $B_{\la+\varep}$ of the form
$\gamma(1)$ for some continuous path $\gamma\dvtx  [0,1]\ra B_{\la+\varep
} $ starting from a point of $B_\la$ and not hitting~$V$; this is an
open set with $V$ as a boundary. Denote also by $W$ another smooth
hypersurface, separating $V$ from $\partial B_{\la} $ and transverse
to $ H_0 $ except on a relative hypersurface. Let now denote by $\OO
'\MMM$ a disjoint copy of the set of past-directed frames
\[
\{(m,\be)\in GL\MMM| \be=(\be_0,\be_1,\dots,\be_d) \mbox{
such that } (m,(-\be_0,\be_1,\dots,\be_d) )\in\OO\MMM\},
\]
and let $O'$, $V'$, ${V}'_0 $ and $W'$ be the subsets of $\OO'\MMM$
corresponding to $O$, $V$, ${V}_0 $ and~$W$. The equivalence relation
\[
(m,(\be_0,\be_1,\dots,\be_d) )\in V \sim (m,(-\be_0,\be_1,\dots
,\be_d) )\in V'
\]
defines a manifold structure on the quotient space $(O\cup{V})\sqcup
(O'\cup{V}') /\sim$, which we denote by $ \mcE$.  Note that $ \mcE
$ is compact and that its volume is in between $2 \Vol(B_\la)$ and $2
\Vol(B_{\la+\varep})$. Write $ \mcV$ for the image in $ \mcE$ of $
V$, and $ \mcV_0 $ for the image in $ \mcE$ of $ V_0 $; define the
primed sets $\mcV'$ and $\mcV_0'$ accordingly.

\begin{rem} \label{rem.geodfE}
The geodesic flow is naturally well defined on $\mcE\bs\mcV_0 $,
getting instantly from $ O$ to $O'$ or from $O'$ to $O$ at its
crossings of $\mcV\bs\mcV_0 $. Indeed by the above definition, for
any $\Phi\in\mcV\bs\mcV_0$, either $H_0(\Phi)$ points outwards
seen from $O$ and inwards seen from $O'$, or $H_0(\Phi)$ points
inwards seen from~$O$ and outwards seen from~$O'$. There is, however,
no a priori convenient way to extend the geodesic flow on $\mcV_0$.
This is the reason why we need to take care of this exceptional set.
\end{rem}

We define the \textit{modified relativistic diffusion} on
the compact manifold $ \mcE$ as follows.
Let $a \dvtx  B_{\la+\varep} \ra[0,1]$ be a smooth function equal to $1$
on $B_\lambda$, and whose vanishing set is exactly the closed part $
\mathfrak{C} $ of $ \mcU$ in between $W$ and $V$ [this means that $
\mathfrak{C} $ is the union of the trajectories $(\gamma_s)_{s\in
(0,1)}\subset\mcU$ of continuous paths $ \gamma$ such that $ \gamma
_0\in W$, $ \gamma_1\in V$, and $(\gamma_s)_{s\in(0,1)}$ does not
intersect the oriented hypersurface $W\cup V$ ].  We extend to $ \mcE
$ the restiction of $ a $ to $O\cup{V}$, by setting $ a(\be') = a(\be
)$ for $ \be'= (m,(-\be_0,\be_1,\dots,\be_d) )\in\OO'\MMM$ and
$ \be= (m,(\be_0,\be_1,\dots,\be_d) )\in\OO\MMM$. We define the
generator of the modified diffusion to be the following variant of $
\mcG_{\Theta} $:
%
%e5.3 ###
\begin{equation}
\label{f.G}
\mcG:= H_0 + \frac{1}{2} \sum_{j=1}^d V_j (a \Theta V_j ) .
\end{equation}
Denote by $ \Vol_{\mcE}$ (resp., $ \Vol_V $, $\Vol_W$) the natural
volume element on $\mcE$ (resp., $V$, $W$).

\begin{lem}
\label{lem.modifE}
For $\Vol_{\mcE}$-almost all starting point $ \Phi_0\in\mcE$, the
modified relativistic diffusion is a well-defined $ \mcE$-valued
process having an almost-surely infinite lifetime.
\end{lem}

\begin{pf}
This modified diffusion has generator $ \mcG_{\Theta} $ in $
B_\lambda$ and in its mirror copy $ B'_\lambda$, and reduces to the
geodesic flow in the region $\{a=0\}$ in between~$W$ and $W'$. After
Remark~\ref{rem.geodfE}, we need first make sure that the set $\mcV
_0\cup\mcV_0'$ of bad points is polar.

Let $\mcN$ and $\mcN'$ be the orbits in the region $\{a=0\}$ of $\mcV
_0$ and $\mcV_0'$ by the geodesic flow. They have, as a consequence of
Lemma~\ref{lem.hyps}, null $\Vol_\mcE$-measure. But as the modified
diffusion started from any $\Phi_0\in\{a>0\}$ is hypoelliptic, its
hitting distribution of $W\cup W'$ has a density \wrt$\Vol_{W\cup
W'}$. It follows that the modified diffusion, started from any point of
$\Phi_0\bs(\mcN\cup\mcN')$, will almost surely never hit $\mcN
\cup\mcN'$, proving that this $\mcE$-valued process is well defined.

It can behave in two ways as it approaches its lifetime: either
crossing infinitely many times $\mcV$, or remaining eventually in a
compact subset of $O$ or $O'$. In the latter case, its projection on
$\MMM$ is a (future or past-directed) timelike path confined in a
compact subset of $O$. As such it has a cluster point at which the
strong causality condition cannot hold, preventing $\MMM$ from being
strongly causal, a contradiction.

In the former case, either the path eventually remains in the region $\{
a=0\}$, or it performs before some finite proper time an infinite
number of crossings from $W\cup W'$ to $\mcV$. Since the geodesic flow
does not explode in $\{a=0\}$, we are left with the latter possibility.
It cannot lead to explosion either, since the geodesic flow needs a
traveling time bounded away from $0$ to travel from $W\cup W'$ to $\mcV$.
\end{pf}
Note that the volume measure $\Vol_{\mcE}$ of the compact manifold $
\mcE$ is an invariant finite measure for the modified diffusion.

%%--------------------------------------------------------------%%
%s5.4 ###
\subsection{\texorpdfstring{Crossing times and escape rate of $\Theta$-diffusions}{Crossing times and escape rate of Theta-diffusions}}
\label{SectionCrossingTimes}
%%--------------------------------------------------------------%%

Fix a reference point $\Phi_{\mathrm{ref}}\in\OO\MMM$, and set
$\mcD(\cdot) = \mcD(\Phi_{\mathrm{ref}} , \cdot)$. Let us
emphasize that $\mcD$ is a two-points function, so it is easy to pass
from $\mcD(\Phi_{\mathrm{ref}} ,\Phi)$ to $\mcD(\Phi_0 ,\Phi)$,
or the other way round, using the triangle inequality, for any $\Phi
_0\in\OO\MMM$.

Given an increasing sequence $(R_n)_{n\geq1}$ of positive reals, set
$\tau_0=0$ and associate to each $R_n$
\[
\mbox{the exit time $\tau_n$ from the box } B^{(n)} := \{\mcD\leq
R_n \}.
\]
It takes the diffusion an amount of proper time $(\tau_n-\tau_{n-1})$
to go from the box $B^{(n-1)}$ to the box $B^{(n)}$. The strategy in
\cite{HsuQin} is to estimate $ \PP_{\Phi}(\tau_n-\tau_{n-1}\leq
t_n)$ for\vadjust{\goodbreak} a suitably chosen deterministic sequence $\{t_n\}_{n\geq0}$
of increments of time. Set for $ n\ge1 $
\[
T_n := \sum_{k=1}^n t_k  \quad \mbox{and}\quad  r_n := R_n-R_{n-1}.
\]
If one can show that
%
%e5.4 ###
\begin{equation}
\label{BorelCantelli}
\sum_{n\geq1} \PP_{\Phi}(\tau_n-\tau_{n-1}\leq t_n) <\infty
\end{equation}
for a convenient choice of the sequences $(R_n)_{n\geq1}$ and
$(T_n)_{n\geq1}$, then the Borel--Cantelli lemma tells us that the
diffusion does not exit $B^{(n)}$ before time $ T_n $, for $n$ large
enough, preventing explosion. Following~\cite{HsuQin}, we are going to
consider the events
\[
E_n := \{\tau_n-\tau_{n-1}\leq t_n , \tau_n\leq T_n\},
\]
so as to be able to use our modified process run backwards from the
\textit{fixed} time $T_n $, when estimating the probability that the
process crosses from~$B^{(n-1)}$ to $B^{(n)}$ not too fast. Lemma $2.1$
of~\cite{HsuQin} (an application of the Borel--Cantelli lemma)
justifies that considering these events leads to the same nonexplosion
conclusion as~\eqref{BorelCantelli}. We recall it here for the
reader's convenience.

\begin{lem} [(\cite{HsuQin})]
\label{LemmeFeinte}
Fix $ \Phi\in\OO\MMM$. If $ \sum_{n\geq1} \PP_{\Phi
}(E_n)<\infty$, then there exists $\PP_{\Phi}$-almost-surely $ \delta$ \st$
\tau_n\geq T_n-\delta,$ for all $n\geq1$.
\end{lem}

We shall use the results of Sections~\ref{SubsectionSubRiemannian} and
\ref{SectionModifiedProcess} to prove the fundamental estimate of
Proposition~\ref{FundamentalCrossingEstimate} below. Given any compact
subset $B$ of $ \OO\MMM$, denote by $\PP_B$ the law of the
relativistic diffusion in $\OO\MMM$ started under the uniform
probability in~$B$
\[
\PP_B(\cdot) = \frac{1}{\Vol(B)}\int_B\PP_{\Phi}(\cdot) \Vol
(d\Phi).
\]
Similarly, and given any compact subset $ A $ of $ \mcE$, write $ \QQ
_A $ for the law of the modified $\Theta$-diffusion in $ \mcE$
started under the uniform probability in $A$.

\begin{prop}
\label{FundamentalCrossingEstimate}
There exists a constant $C$ \st we have, for any $n\geq1$\emph{,}
\begin{eqnarray*}
&&\PP_{B^{(1)}} (\tau_n-\tau_{n-1}\leq t_n , \tau_n\leq T_n ) \\
&&\qquad \leq C\frac{\Vol(B^{(n)})}{\Vol(B^{(1)})} \frac{T_n \sqrt{\widehat
{\Theta}_n/t_n}}{(r_n-1-4t_n)} \exp\biggl[{{- \frac{(r_n-1-4t_n)^2}{32
\widehat{\Theta}_n t_n}}} \biggr] ,
\end{eqnarray*}
where $ \widehat{\Theta}_n$ denotes the supremum of $ \Theta$ over
the box $ \{\mcD\leq R_n+1 \}$.
\end{prop}

 The proof mimics Takeda's original proof, as adapted by Hsu and
Qin in~\cite{HsuQin}, with the noticeable difference that we are
working with a nonsymmetric, nonelliptic diffusion.\vadjust{\goodbreak}

\begin{pf*}{Proof of Proposition~\ref{FundamentalCrossingEstimate}}
We start by embedding the box $B^{(n)}$ into the set $\mcE^{(n)}$
constructed in Section~\ref{SectionModifiedProcess}, with $ \la=R_n $
and $ \varep= \frac{1}{2} $, say. From now on we
work on the path space over $\mcE^{(n)}$ and use the coordinate
process~$X$, whose filtration is denoted by $(\mcF_s)_{s\geq0}$. We
still denote by $\tau_n$ the exit time from (the image in $\mcE
^{(n)}$ of) $B^{(n)} $; the event
\[
E_n := \{\tau_n-\tau_{n-1}\leq t_n, \tau_n\leq T_n\}
\]
belongs to $\mcF_{\tau_n}$. As explained above in Section \ref
{SectionMethod}, the proof has two main ingredients, the first of which
is inequality~\eqref{FundamentalBasicInequality} below, where $\QQ
_{\mcE^{(n)}}$ denotes the distribution of the modified $\Theta
$-diffusion in $\mcE^{(n)}$, with generator $\mcG$ given in~\eqref{f.G}.

As the $\Theta$-diffusion and the modified $\Theta$-diffusion have
the same law before the stopping time $\tau_n$, we have $\PP
_{B^{(n)}}(E_n)=\QQ_{B^{(n)}}(E_n) \leq2 \QQ_{\mcE^{(n)}}(E_n)$,
and so
%
%e5.5 ###
\begin{equation}
\label{FundamentalBasicInequality}
\PP_{B^{(1)}}(E_n) \leq 2 \frac{\Vol(B^{(n)})}{\Vol(B^{(1)})} \QQ
_{\mcE^{(n)}}(E_n) ,
\end{equation}
by the obvious inequality $ \PP_{B^{(1)}}(E_n)\leq\frac{\Vol
(B^{(n)})}{\Vol(B^{(1)})} \PP_{B^{(n)}}(E_n)$. The second ingredient
involves the Lyons--Zheng decomposition of $ \mcD(X_s)$ under $\QQ
_{\mcE^{(n)}}$. As~$\mcD$ is not a priori sufficiently regular to use
It\^{o}'s formula, we apply it to its smooth approximation $F$
constructed in Proposition~\ref{CorApproxTimeFunction} (with $R=R_n$
and $\eta=\frac{1}{2}$). As the process
$(X_{T_n-s})_{0\leq s\leq T_n}$ is under $\QQ_{\mcE^{(n)}}$ a
homogeneous diffusion process with generator $\mcG^* = -H_0 + \frac
{1}{2} \sum_{j=1}^d V_j (a \Theta V_j )$, it follows from It\^{o}'s
formula that there exist two martingales $(M_s)_{0\leq s\leq T_n}$ and
$(\widetilde{M}_s)_{0\leq s\leq T_n}$, \wrt the forward and backward
filtrations of the process, respectively, \st
\begin{eqnarray*}
F(X_s) &=& F(X_0) + M_s + \int_0^s \mcG F(X_r) \,dr, \\
F(X_s) &=& F\bigl(X_{T_n-(T_n-s)}\bigr) = F(X_{T_n}) + \widetilde{M}_{T_n-s} +
\int_s^{T_n} \mcG^*F(X_r) \,dr,
\end{eqnarray*}
with
\begin{eqnarray}
\label{ControlBracket}
\langle M\rangle_s &=& \sum_{j=1}^d \int_0^s a(X_r) \Theta(X_r)
|V_jF |^2(X_r) \,dr \leq4 \widehat{\Theta}_n s, \nonumber\\[-8pt]\\[-8pt]
\langle\widetilde{M} \rangle_s &=& \sum_{j=1}^d \int
_0^sa(X_{T_n-r}) \Theta(X_{T_n-r}) |V_jF |^2(X_{T_n-r}) \,dr \leq4
\widehat{\Theta}_n s.\nonumber
\end{eqnarray}
Setting $M'_s:= \widetilde{M}_{T_n-s} $ and noting that $\mcG-\mcG^*
= 2H_0$, we thus have
%
%e5.6 ###
\begin{equation} \label{f.F(X)MM}
d (F(X_s) ) = d \biggl(\frac{M_s+M'_s}{2} \biggr) + H_0F(X_s) \,ds
\end{equation}
with a controlled drift term $ |H_0F|\leq2 $, by Proposition \ref
{CorApproxTimeFunction}. By construction, we have
\[  %\label{f.F(X)}
\sup_{0\leq s\leq t_n} |F(X_{\tau_{n-1}+s})-F(X_{\tau
_{n-1}}) | \geq r_n -1
\]
on the event $ E_n $, where $X$ hits the set \mbox{$\{F\geq R_n-\frac
{1}{2}\}$} in the time interval $[\tau_{n-1}, \tau_{n-1}+t_n]$. To
control the $\QQ_{\mcE^{(n)}}$-probability of $E_n$, we use Hsu and
Qin's trick. Cut the interval $[0,T_n] = \bigcup_{k=1}^{\ell_n}
[(k-1)t_n,k t_n ] $ into $ \ell_n := T_n/t_n $ sub-intervals of length
$ t_n$ (to lighten the notations, we shall neglect the fact that $ \ell
_n $ may not be an integer$ $; this fact causes no trouble but
notational), and write on each event $ \{(k-1)t_n\le\tau_{n-1}\le k
t_n \}$
\[
F(X_{\tau_{n-1}+s})- F(X_{\tau_{n-1}}) = F(X_{\tau_{n-1}+s})-
F(X_{kt_n})+ F(X_{kt_n})- F(X_{\tau_{n-1}}).
\]
This simple remark shows that the event $ \{\sup_{0\leq s\leq
t_n} |F(X_{\tau_{n-1}+s}) -F(X_{\tau_{n-1}}) |\geq r_n -1 \}$
is included in one of the $ \ell_n $ events $  \{\sup_{0\leq
|s|\leq t_n} |F(X_{kt_n+s}) -F(X_{kt_n}) |\geq \frac
{r_n-1}{2} \}$, where $1\leq k\leq\ell_n$. By~\eqref{f.F(X)MM} and
the inequality $|H_0F|\leq2$, the $k$th of these events
is included in the union $ A_k\cup\widetilde A_k $, where
\[
A_k := \biggl\{\sup_{0\leq|s|\leq t_n} |M_{kt_n+s}-M_{kt_n}
|\geq\frac{r_n-1}{2} - 2t_n \biggr\}
\]
and
\[
\widetilde A_k := \biggl\{\sup_{0\leq|s|\leq t_n} |\widetilde
{M}'_{kt_n+s}-\widetilde{M}'_{kt_n} |\geq\frac{r_n-1}{2} -
2t_n \biggr\}.
\]
Let $W$ be a Brownian motion  defined on some probability space $(\Omega,\mcF,\PP
)$. By~\eqref{ControlBracket} we have
\begin{eqnarray*}
\QQ_{\mcE^{(n)}}(A_k) &\leq& 2 \PP\biggl(\sup_{0\leq s\leq t_n}|W_s|\geq\frac{r_n-1-4t_n}{4\sqrt{\widehat{\Theta}_n}}
\biggr)\\
&\leq& \frac{C \sqrt{\widehat{\Theta}_n/t_n}}{r_n-1-4t_n} \exp
\biggl(-\frac{(r_n-1-4t_n)^2}{32 \widehat{\Theta}_n t_n} \biggr)
\end{eqnarray*}
for some positive constant $ C $; the same identity holds for
$\widetilde A_k$, using~\eqref{ControlBracket}. Summing over $ k $ and
using inequality~\eqref{FundamentalBasicInequality} yields the
statement of the proposition since $E_n\subset\bigcup_{k=1}^{\ell
_n}(A_k\cup\widetilde{A}_k)$.
\end{pf*}

This key proposition being proved, it becomes easy to prove Theorem
\ref{MainThm}.

\begin{pf*}{Proof of Theorem~\ref{MainThm}}
Taking $R_n=2^{n+5}$ and $ t_n\leq2^{n+1}$ in Proposition \ref
{FundamentalCrossingEstimate}, so that $T_n\leq2^{n+2}$, we get for
any $n\geq1 $
%
%e5.7 ###
\begin{eqnarray}
\label{EquationEstimate}
\PP_{B^{(1)}}(E_n) &=& \PP_{B^{(1)}} (\tau_n-\tau_{n-1}\leq t_n, \tau
_n\leq T_n ) \nonumber\\[-8pt]\\[-12pt]
&\leq& C \frac{\Vol(B^{(n)})}{\Vol(B^{(1)})} \sqrt{\frac
{\widehat{\Theta}_n}{t_n}} \exp\biggl[{- \frac{4^n}{\widehat{\Theta}_n
t_n}} \biggr] .\nonumber
\end{eqnarray}
Specifying the choice of $ t_n $ by setting
\[
t_n := \min\biggl\{{2^{n+1}} , \frac{4^{n-1}}{ ( 1+ \log^+ [ \widehat
{\Theta}_n \Vol(B^{(n)}) ] ) \widehat{\Theta}_n} \biggr\},
\]
the right-hand side of~\eqref{EquationEstimate} is seen to be bounded
above by a constant multiple of $2^{-n}$, ensuring as a consequence the
convergence of the series $ \sum_{n\geq1} \PP_{B^{(1)}}(E_n)$.
Indeed, we get from~\eqref{EquationEstimate}, with the above $t_n$,
\[
\PP_{B^{(1)}}(E_n) \le C' \Vol\bigl(B^{(n)}\bigr) \sqrt{\frac{\widehat
{\Theta}_n^2 \log[ \widehat{\Theta}_n \Vol(B^{(n)}) ]}{4^n}} e^{-
{4}\log[ \widehat{\Theta}_n \Vvol(B^{(n)}) ]} \le C''/2^n.
\]
[Ignoring the trivial case $\Theta\equiv0$, we can suppose without
loss of generality that we have $\widehat{\Theta}_n \Vol(B^{(n)})
\geq3$ for $n$ large enough.] Note that the above choice of time
increments $t_n$ is simpler than Hsu and Qin's choice in \cite
{HsuQin}; there is in particular no need to introduce their auxiliary
function $ h(R) \equiv\log\log R $, to get Grigor'yan's criterion, if
the second upper bound of their Section 3 is not used.

To conclude that the $\Theta$-diffusion does not explode we need to
check that $ T_n = \sum_{k=1}^n t_k$ increases to infinity. For the
above choice of time increments~$t_n$, we have $\PP_{B^{(1)}}$-almost-surely,
for $n$ larger than some $n_0$, and for a positive universal constant~$c$,
\begin{eqnarray}\label{f.minTn}
T_n& \ge&\sum_{k= n_0}^n \min\biggl\{2^{k+1} , \frac{4^{k-1}}{\Theta
_{2^{k+5}+1} (\log^+ [ \Theta_{2^{k+5}+1}\Vol(B_{2^{k+5}}) ] + 1 )}\biggr\}
\nonumber\\[-9pt]\\[-9pt]
&\ge &c \int_{2^{n_0+1}}^{2^n} \min\biggl\{8, \frac
{r}{\Theta_r \log[\Theta_r \Vol(B_r) ] } \biggr\} \,dr .\nonumber
\end{eqnarray}
Leaving aside the trivial case $\Theta\equiv0$ and recalling that the
map $r\mapsto\Theta_r = \max_{B_r} \Theta$ is
nondecreasing, we can suppose \wlg that $\Theta_r\geq3$. The
divergence of the sequence $(T_n)$ is then granted by the integral criterion
\[
\int^\infty\min\biggl\{8, \frac{r}{\Theta_r \log
[\Theta_r \Vol(B_r) ] } \biggr\} \,dr = \infty.
\]
As $\Theta_r$ increases, this condition is equivalent to
\[
\sum_{n\geq1} \min\biggl\{8, \frac{n}{\Theta_n
\log[\Theta_n \Vol(B_n) ] } \biggr\} = \infty,
\]
that is to
\[
\sum_{n\geq1} \frac{n}{\Theta_n \log[\Theta_n \Vol(B_n) ] } =
\infty,\vadjust{\goodbreak}
\]
since the former holds obviously if an infinite number of terms were
larger than 8. The previous condition is equivalent to condition~(\ref{ConditionMainThm}) of Theorem~\ref{MainThm}.

Using the Borel--Cantelli lemma under the form of Lemma \ref
{LemmeFeinte}, it follows that we have
%
%e5.8 ###
\begin{equation}
\label{f.result}\quad
\PP_{B^{(1)}} \Bigl( \sup_{0\leq s\leq T_n-\delta} \mcD(\Phi
_s)\leq2^{n+5} \mbox{ for any large enough } n \Bigr) = 1 ,
\end{equation}
so $\sup_{0\leq s\leq t} \mcD(\Phi_s)<\infty$, for all
$t>0$, since $T_n$ increases to $\infty$. Would a realization of the
path $\Phi_s$ explode by time $t$, its projection in $\MMM$ would
provide a timelike path with an accumulation point [for it stays in the
projection of a compact set by Hypothesis \hyperlink{hypH}{(\textup{H})}], contradicting
the strong causality assumption on $\MMM$.

To prove that the same happens under any $\PP_{\Phi_0}$, notice that
since the nonexplosion event $E$ belongs to the invariant $\sigma
$-algebra, the function $\OO\MMM\ni\Phi\mapsto\PP_\Phi(E)$ is
$\mcG_\Theta$-harmonic, hence continuous, as $\mcG_\Theta$ is
hypoelliptic. It follows that since
\[
\PP_{B^{(1)}}(E) = \frac{1}{\Vol(B^{(1)} )}\int_{B^{(1)}} \PP
_{\Phi}(E) \Vol(d\Phi) ,
\]
the probability $\PP_{\Phi}(E)$ must be equal to $1$ for all $\Phi
\in B^{(1)}$. But as the ball~$B^{(1)}$ was arbitrarily chosen, $\PP
_{\Phi}(E)$ is identically equal to $1$ everywhere.
\end{pf*}

%------------------------------%
%s5.5 ###
\subsection{Upper rate function}
\label{SectionUpperRateFunction}
%------------------------------%

Using essentially the same reasoning as in Section~4 of~\cite{HsuQin},
the above proof yields almost for free the upper rate function for the
$\Theta$-diffusion given in Corollary~\ref{CorMainThm}. See also
\cite{GrigoryanEscapeRate} for related results. We keep the preceding
notation.

\begin{pf*}{Proof of Corollary~\ref{CorMainThm}}
We follow the argument of~\cite{HsuQin}, Section 4, making sure that
it works here as well with our choice for $t_n$, and without their
auxiliary function $\log\log$. Suppose first $\Theta$
nonidentically null and recall inequality~\eqref{f.minTn}, in which
we can forget to take the minimum with 8, by Proposition \ref
{PreparatoryLemma} below. By~\eqref{f.result}, this yields,
almost-surely, the inequality
\[
\sup_{0\leq s\leq c h^{-1}(2^n)-\delta} \mcD(\Phi
_s)\leq2^{n+5},
\]
that is
\[
\sup_{0\leq s\leq c h^{-1}(R)-\delta} \mcD(\Phi_s)\leq
32 R,
\]
for large enough $ R $. Letting $ R=h ( (t+\delta)/c )$, this entails
  $\sup_{0\leq s\leq t} \mcD(\Phi_s)\leq32 h (
(t+\delta)/c )$,  hence  $\sup_{0\leq s\leq t} \mcD
(\Phi_s)\leq32 h(C t)$,  for large enough $ t $.\vspace*{1.5pt} This shows the
claim under the probability $\PP_{B^{(1)}}$, and then under $ \PP
_{\Phi_0}$ as well, by the same argument already used at the end of
the proof of Theorem~\ref{MainThm}. Finally, in the geodesic case $
(\Theta\equiv0 )$, the same holds with $ T_n\ge c 2^n = c h(2^n)$.
\end{pf*}

%%-------------------------------------------------------------------------------%
%s5.6 ###
\subsection{Estimates of the volume of the sub-Riemannian boxes and
application}
\label{sec.Thex}
%%-------------------------------------------------------------------------------%

Let us begin with a crude lower estimate of the volume of the
boxes~$B_r$ based on the vertical expansion in the $\mathit{SO}_0(1,d)$-fiber of $\OO
\MMM$, without taking into account the horizontal expansion which
depends on the curvature of the base Lorentzian manifold $\MMM$. We
used this lower bound in the proof of Corollary~\ref{CorMainThm}.
\begin{prop} \label{PreparatoryLemma}
We have $ \liminf_{r\to\infty} \frac{\log{\Vvol}(B_r)}{r} \ge d-1 $.
\end{prop}
\begin{pf}
Fix a relatively compact neighborhood $\mcU$ of $m_0$ in $\MMM$,
above which $\OO\mcU$ is trivialized in $\mcU\times \mathit{SO}_0(1,d)$.
Assume \wlg that $\Phi_0$ corresponds to $(m_0,\mathbf{1})$. By the
ball-box theorem (see, e.g.,~\cite{Montgomery}), the box $B_r = \{\mcD
\le r\}$ contains a neighborhood $ \mcV\times B(\mathbf{1},\varep)$ of
$\Phi_0 $, for some $\varep>0$ and for $r$ larger than some fixed
$r_1$. Using this argument a finite number of times, together with the
triangle inequality for $\mcD$, we see that the box $\{\mcD\le r\}$
contains any neighborhood $\mcU\times B(\mathbf{1},\rr)$ of $\Phi_0$,
for any $\rr>0$, provided $r$ is large enough, say no less than
$r_0=r_0(\mcU,\rr)$. Take $\rr$ larger than the diameter of $
\mathit{SO}(d)$.

We easily see that the boxes $\{\mcD\le r\}$ dilate in the vertical
directions $V_1,\ldots,V_d$ with speed $r$, as $r$ increases. So $\{
\mcD\le r\}$ contains the product of~$\mcU$ by the ball of radius
$(r-r_0)$ in $\mathit{SO}_0(1,d)$ for $r$ large enough. This provides a lower
bound on ${\Vol} (\{\mcD\le r\} )$ by some constant multiple of the
volume of the hyperbolic ball of radius $(r-r_0)$, from which it
follows that there exists some positive constant $c$ \st$\log{\Vol
}(B_r) \ge(d-1) r+\log c $, for $r$ large enough.
\end{pf}

To close this work, we give a nonexplosion criterion involving only the
geometry of $\MMM$, rather than the geometry of $\OO\MMM$ as it
appears in Theorem~\ref{MainThm} through the sub-Riemannian boxes $B_r$.

\begin{prop}
\label{pro.estT}
Fix $\Phi_0=(m_0,\be_0)\in\OO\MMM$, and define the \textup
{$S_{\Phi_0}$-radius} $ \rho^S_{\Phi_0}(m)$ of any $ m\in
\MMM$ as the infimum of the $S_{\Phi_0}$-length of $C^1$ paths
joining $m_0$ to $m$. Define the \textup{$S_{\Phi_0}$-ball}
$ B^S_{\Phi_0}(r)$ of radius $ r $ as the set $B^S_{\Phi_0}(r):= \{
m\in\MMM | \rho^S_{\Phi_0}(m)\leq r \}$, and set
\[
V^S(r) := \Vol_{\MMM} ( B^S_{\Phi_0}(r) ).
\]
Then there exists a constant $C$ such that we have for all $r>0$
\[
\log\Vol(B_r) \le C + (d-1) r + \log V^S(C e^r).
\]
\end{prop}

Note that the $S_{\Phi_0}$-balls $B^S_{\Phi_0}(r)$ and their volume
depend only on the choice of $\Phi_0=(m_0,\be_0)\in\OO\MMM$ and on
the geometry of $\MMM$. We noticed indeed in Section \ref
{SectionbCompleteness} that the $S_{\be_0}$-length of a path in $\MMM
$ started from $m_0$ is the Euclidean length of its anti-development in
$ (T_{m_0}\MMM,\be_0 )$.\vadjust{\goodbreak}

\begin{pf*}{Proof of Proposition~\ref{pro.estT}}
By the definitions in Sections~\ref{SectionbCompleteness} and \ref
{SubsectionSubRiemannian}, the b-distance of $\Phi_0$ to any $ \Phi
\in\OO\MMM$ is not larger than $\mcD_{\Phi_0}(\Phi)$, so
$B_r\subset B^b(\Phi_0 ; r)$, where $B^b$ denotes the ball in $\OO
\MMM$ of the b-metric. Vertically, that is to say in the frame $ \tau
_{0\to s}^\gamma(\Phi_0)$ parallely transported along a minimizing
curve~$\gamma$, the maximal hyperbolic distance reached by the
velocity component~$\dot m_s$ of~$\gamma_s$ is $s$, which is
responsible for a maximal vertical volume $ \mcO(e^{(d-1) r})$.

Having accelerated till reaching a maximal velocity $\mcO(e^r)$, a
minimizing curve in $B^b(\Phi_0 ; r)$ can perform a maximal horizontal
displacement $\mcO(e^r)$. Hence we have the inclusions
\[
B^S_{\Phi_0}(r)\subset\pi_0 (B^b(\Phi_0 ; r) ) \subset B^S_{\Phi
_0} (\mcO(e^r) ),
\]
and so $\Vol(B_r) \le C e^{(d-1) r} V^S(C e^r)$.
\end{pf*}

Applying Proposition~\ref{pro.estT} to the integral condition of
Theorem~\ref{MainThm} yields in the case of a bounded $\Theta$ the
nonexplosion criterion $\int^\infty\frac{r \,dr}{ r + \log
V^S(e^r)} =\infty$. Using the increasing character of the map $
(r\mapsto V^S(e^r) )$, discretizing and distinguishing whether or not
there are infinitely many $n$ such that $\log V^S(e^n)\leq n$, we
easily see that this condition is equivalent to the condition $\int
^\infty\frac{r \,dr}{\log V^S(e^r)} =\infty$.

\begin{cor}
\label{cor.estT}
Let $(\MMM,g)$ be a strongly causal Lorentz manifold satisfying the
Completeness hypothesis \hyperlink{hypH}{(\textup{H})} and the volume growth
condition $ \int^\infty\frac{r \,dr}{\log V^S(e^r)} = \infty$.
Then all $\Theta$-diffusions with a bounded $\Theta$ are
stochastically complete.
\end{cor}

It is easy to see that this volume growth integral criterion does not
depend on the choice of $\Phi_0\in\OO\MMM$. Contrary to Proposition
\ref{PreparatoryLemma}, it relies on the horizontal expansion and not
on the vertical expansion. This criterion does not apply to G\"{o}del's
universe, for which $\log V^S(e^r)$ is of order $e^r$; the nonexplosion
criterion of Section~\ref{SectionPerfectFluids} covers the case of
that spacetime. Corollary~\ref{cor.estT} applies, for example, to
Lorentz manifolds which are topologically $\RR^{1+d}$ and have a
pseudo-metric $g$ such that $g, g\1$ and the first-order derivatives of
$g$ \wrt the canonical coordinates are bounded, since then $\log
V^S(e^r)$ is of order $r$, as is the case in Minkowski spacetime.

\section*{Acknowledgments}

We thank E. Tr\'{e}lat for his guidance in the realm of control theory
and A. Oancea and P. Pansu for their help in proving Lemma~\ref{lem.hyps}.

%suskaldyti doi

% imsref loaded by dianan, 2011-05-13 10:50:49
%

\printaddresses

\end{document}